\documentclass[a4paper, 11pt, abstract=true, oneside]{scrartcl}  
\usepackage[authoryear,sort&compress]{natbib}

\usepackage[utf8x]{inputenc}
\usepackage[english]{babel}
\usepackage{graphicx, amsfonts,amsmath,amssymb,amsthm,bm,mathrsfs}
\usepackage{commands}
\usepackage{longtable,fancybox}
\usepackage{comment,geometry}
\usepackage{booktabs, multirow}
\usepackage{makecell}
\usepackage{threeparttable}
\usepackage{caption, subcaption}
\usepackage{enumitem}

\usepackage{pgfplots}
\usepackage{tikz}
\usetikzlibrary{positioning}
\usetikzlibrary{shapes.geometric}
\usetikzlibrary{fit}
\usetikzlibrary{external}
\def\Halmos{\mbox{\quad$\square$}}

\usepackage{algorithm,algpseudocode}

\AtBeginDocument{
	\newgeometry{
		textheight=9in,
		textwidth=6.5in,
		top=1in,
		headheight=16.49675pt,
		headsep=25pt,
		footskip=30pt
	}
}
\author{Mathis Brichet, Axel Parmentier, Maximilian Schiffer}

\title{\LARGE \textbf Optimizing a Worldwide-Scale Shipper Transportation Planning in a Carmaker Inbound Supply Chain}

\date{October 12, 2026}

\newcommand{\note}[1]{\caption*{\small \textit{Notes} #1}}

\begin{document}

\maketitle

\begin{abstract}

We study the shipper-side design of large-scale inbound transportation networks, motivated by the global supply chain of the carmaker Renault. 
We formalize the \emph{Shipper Transportation Planning Problem} (STPP), which integrates discrete flow consolidation via explicit bin-packing, time-expanded routing, and operational regularity constraints. 
To solve this high-complexity combinatorial problem at an industrial scale, we propose a tailored \emph{Iterated Local Search} (ILS) metaheuristic. 
The algorithm combines large-neighborhood search with MILP-based perturbations and leverages bundle-specific decompositions to obtain effective search guidance. 
Computational experiments on real industrial data involving more than 700,000 commodities and 1.2 million arcs demonstrate that the ILS achieves an average gap of 7.9\% to the best available lower bound. 
The results reveal a 23.2\% cost-reduction potential compared to legacy planning benchmarks. 
Most significantly, the proposed framework is currently deployed in production at Renault, where it supports weekly strategic decisions.
Internal evaluations indicates realizable cost savings estimated at approximately €20 million per year. 
Our analysis yields key managerial insights: we demonstrate that explicit 1D bin-packing is a critical step forward for realistic consolidation modeling, that transport regularity offers a robust balance between cost and stability, and that high-volume global networks benefit significantly from in-house strategic planning over third-party outsourcing. 
To the best of our knowledge, this is the first work to successfully solve a shipper-side transportation design problem at this magnitude.

\end{abstract}

\section{Introduction}\label{sec:Introduction}

Designing an efficient inbound supply chain is a critical task for manufacturing firms that operate at massive scale. 
The inbound supply chain connects a company’s suppliers to its production sites, and its configuration must often be revised in response to product introductions, demand shifts or the need to mitigate disruptions.
In industries with geographically distributed suppliers and industrial sites, this task becomes especially complex.
While many firms outsource logistics operations to third-party logistics providers (3PLs) to leverage shared transportation assets and reduce operational complexity, this outsourcing approach is not universally applicable nor beneficial.
Some manufacturers operate supply chains of such magnitude and intricacy that they prefer to retain full control over their logistics planning.
This is particularly true for large automotive manufacturers, who can achieve economies of scale internally.
For these firms, outsourcing would not only reduce flexibility and transparency but also forgo substantial cost savings and introduce operational risk.

This paper arises from a close industrial collaboration with Renault and addresses a real-world strategic problem: the long-term planning of the company’s inbound transportation network.
As a global car manufacturer, Renault moves millions of parts across continents each year.
This results in an inbound logistics network of exceptional size and complexity—comprising thousands of suppliers, hundreds of logistics platforms, and dozens of industrial sites—requiring specialized and scalable optimization methods to ensure efficient operations.
In this context, Renault initiated a strategic effort to modernize and strengthen its inbound logistics planning capabilities. 
Recent initiatives have focused on redesigning planning processes around systematic use of operational data and advanced decision‑science tools, resulting in the development of multiple tools adapted to each planning level: strategic, tactical and operational.
Due to the complexity of consolidation mechanisms and the scale of the network, 
Renault required a new generation of planning tools—beyond off‑the‑shelf software or generic outsourcing models—built on bespoke, high‑performance optimization methods capable of handling the full global problem, capturing the system’s specific structural characteristics, incorporating finer consolidation models, and producing high‑quality solutions within reasonable computation times.

In this work, our objective is to produce a provisional transportation plan for parts to be assembled in the next six months.
Each commodity must be routed through Renault’s logistics network in a cost-efficient manner.
The decision involves choosing whether to ship parts directly or to consolidate them through intermediate platforms.
The model is not intended to generate the exact transportation plan ultimately executed over the six-month horizon. 
Rather, its purpose is to identify major consolidation corridors and long-term recurring routing patterns that support carrier procurement contracts, capacity negotiations, and global network evolution decisions. 
Routing and loading decisions are modeled because they determine the effective utilization of transportation services and therefore the economic attractiveness of alternative corridors. 
These decisions serve as planning proxies rather than operational commitments, while detailed shipment scheduling and day-to-day execution remain the responsibility of downstream tactical and operational planning processes.
To allow for optimized decision-making in this context, we formalize the underlying Shipper Transportation Planning Problem (STPP), tailored to the needs of a high-volume industrial shipper, and develop a scalable algorithmic framework, capable of generating high-quality plans on real data from one of Europe’s largest manufacturers.

\subsection{State of the art}\label{subsec:SOTA}

The design and optimization of transportation networks have been long-standing challenges in operations research, especially in the context of logistics and supply chain management.
The foundational modeling framework is that of multicommodity network flow problems, extensively studied since the 1960s \citep{ahuja1988network}.
Early research already recognized the critical importance of flow consolidation in freight transportation \citep{powell1983load}.
The service network design (SND) literature, a broad class of planning problems that captures consolidation-based freight operations, is the appropriate reference point for our problem.
It was formally established in \citet{crainic2000service} and further analyzed by \citet{andersen2009service} and \citet{crainic2020exact}. 
More recently, \citet{crainicNetworkDesignApplications2021a} \citet{crainic202550} provide comprehensive surveys of the methodological landscape.

Within this family of problems, our work extends a scheduled multi-layer service network design problem with unsplittable flows by incorporating explicit bin-packing consolidation, flexible delivery times, additional regularity constraints, and a strategic shipper-oriented viewpoint, all within an extremely large-scale industrial setting. 
Among the various SND variants, we focus primarily on load plan design problems (LPDPs), which share many structural similarities with our formulation while remaining more tractable than full-fledged SND models.

Table~\ref{tab:state_of_the_art} compares our model with representative works in the SND and LPDP literature: [1] \citet{jarrah2009}, [2] \citet{lindsey2016improved}, [3] \citet{eomrecursive}, [4] \citet{he2022exact} and [5] \citet{zhu2014scheduled}.
\begin{table}[t]
    \centering
    \caption{Similarities and differences with the most relevant related work.}\label{tab:state_of_the_art}
    \renewcommand{\arraystretch}{1.3} 
    \resizebox{\textwidth}{!}{%
    \begin{tabular}{p{0.3\textwidth} p{0.175\textwidth} p{0.175\textwidth} p{0.175\textwidth} p{0.175\textwidth} p{0.175\textwidth} p{0.175\textwidth}}
            \toprule
            \textbf{Characteristics} & \textbf{[1]} & \textbf{[2]} & \textbf{[3]} & \textbf{[4]} & \textbf{[5]} & \textbf{Our work} \\
            \midrule
            Method & Heuristic & Heuristic & Heuristic & Exact & Heuristic & Heuristic \\
            Flow type & U-MCF & U-MCF & U-MCF & U-MCF & U-MCF & U-MCF \\
            Unrestricted path set & \texttimes & \texttimes & \checkmark & \checkmark & \texttimes & \checkmark \\
            Network nodes (time-space) & 5100 & 5400 & 1000 & 4300 & 210 & 63 500 \\
            Network arcs (time-space) & 210 000 & 864 000 & 248 000 & 135 000 & 6100 & 1 190 000 \\
            Number of commodities & 680 & 60 000 & 40 000 & 1000 & 1000 & 700 000 \\
            Time horizon & 1 week & 4 days & 1 day & 3 hours & 1 week & 6 months \\
            Bin-Packing consolidation & \texttimes & \texttimes & \texttimes & \texttimes & $\sim$ & \checkmark \\
            Regularity type & T-G-IT & T-G-IT & \texttimes & \texttimes & G & T-G \\
            Resource management & \checkmark & \checkmark & \texttimes & \texttimes & \texttimes & \texttimes \\
            Flexible delivery time & \texttimes & \texttimes & \texttimes & \texttimes & \texttimes & \checkmark \\
            Point of view & Carrier & Carrier & Carrier & Carrier & Carrier & Shipper \\
            Decision scope & Tactic & Tactic & Tactic & Operations & Tactical & Strategic \\
            Lower Bound & \texttimes & \texttimes & \checkmark & \checkmark & \checkmark & \checkmark \\
            \bottomrule
        \end{tabular}
    }
    \vspace{-0.25cm}
    \note{Usage: U = Unsplittable, MCF = Multi-Commodity Flow, T = Time regularity, G = Group regularity, IT~=~In-Tree regularity.}
\end{table}
Table~\ref{tab:state_of_the_art} highlights several structural differences induced by a shipper-oriented perspective relative to the carrier-centric formulations that dominate the SND and LPDP literature. 
Unlike carrier-side models, transportation resources are procured rather than managed, shifting the focus from resource scheduling to transportation-policy design over longer planning horizons. 
The explicit representation of suppliers, destinations, and shipment streams also leads to substantially larger and sparser networks, while direct access to shipment characteristics enables commodity-level consolidation decisions through bin-packing constraints. 
Finally, flexible delivery schedules increase the importance of grouping and temporal-regularity decisions. Collectively, these characteristics define modeling and algorithmic challenges that differ substantially from those addressed in the traditional SND literature.

Exact solution methods for such problems typically rely on mathematical programming techniques, including decomposition methods such as Lagrangian relaxation, Benders decomposition, and Dantzig-Wolfe reformulation \citep{crainic2020exact, boland2017continuous, fontaine2021scheduled, he2022exact} [4].
These methods are effective on moderate-size instances but struggle to scale.
Studies by \citet{gendron2014branch} and \citet{frangioni2017computational} report computation times of several hours even on modest networks, involving a few hundred arcs and commodities.

To address larger systems, the literature increasingly relies on heuristic and matheuristic techniques. 
Metaheuristics such as large neighborhood search or population-based approaches have demonstrated strong empirical performance on realistically sized problems \citep{crainic2021heuristics, kazemzadeh2022node, gendron2018matheuristics, paraskevopoulos2016cycle, jarrah2009} [1]. 
In road transportation planning, for instance, \citet{erera2013improved} and \citet{lindsey2016improved} [2] solve instances involving tens of thousands of arcs and commodities by combining local search strategies with sophisticated descent-based matheuristics. 
More recently, \citet{eomrecursive} [3] proposed a recursive partitioning and batching strategy to improve scalability through problem decomposition.
Our present work, with massive commodity volumes, demands algorithmic solutions that scale even further.

Explicit consolidation decisions have received far less attention. 
To the best of our knowledge, no LPDP formulation integrates bin-packing constraints to determine the number of services required on each arc \citep{bakir2021motor}. 
A partial exception occurs in train blocking and scheduling models \citep{zhu2014scheduled, kienzle2025intermodal}, where commodities are first grouped into blocks, which are then assigned to trains. 
While this two-stage consolidation mimics bin-packing, block structures are precomputed and therefore eliminate flexibility in commodity-to-block assignments. 
The absence of explicit bin-packing constraints is noteworthy: introducing them dramatically increases problem size, complicates decomposition and matheuristic strategies and changes the solution obtained. 
By explicitly modeling discrete 1D bin-packing, we move beyond the standard volume relaxations prevalent in the literature.

Finally, existing LPDP models overwhelmingly adopt a carrier-side viewpoint, embedding resource management constraints and fixed departure schedules \citep{lindsey2016improved}. 
In contrast, a shipper-oriented setting focuses primarily on transportation procurement and the design of consolidation corridors. 
Transportation resources are acquired through a portfolio of carriers rather than operated directly, shifting the planning challenge from resource scheduling to estimating transportation volumes, negotiating capacity requirements, and coordinating freight flows across the network. 
Direct access to shipment-level information also supports longer planning horizons, larger planning scopes, and greater flexibility in delivery timing.
From an algorithmic perspective, these characteristics substantially enlarge the feasible solution space and create stronger interactions between commodities through consolidation effects, thereby alleviating the effectiveness of approaches that rely on enumerating feasible paths. 
They instead motivate solution methods capable of jointly handling large-scale routing, consolidation, and coordination decisions.

Overall, this review highlights a significant gap: while LPDPs and SND problems have been widely studied from the perspective of carriers and third‑party logistics providers, the shipper’s perspective—combined with explicit consolidation decisions and massive commodity volumes for which current approaches remain insufficient—remains largely unexplored.

\subsection{Contributions}

This paper studies the shipper-side design of large-scale inbound transportation networks, focusing on the strategic planning problem faced by Renault. 
Building on the gaps identified in 
literature, we make four main contributions:

\begin{description}
    \item[\textbf{(i) A novel shipper problem setting within the SND family.}] We formalize the \emph{Shipper Transportation Planning Problem} (STPP), which integrates several structural features rarely combined in prior studies: explicit bin-packing consolidation, flexible delivery-time windows, regularity constraints, and massive-scale unsplittable flows. 
    Unlike carrier-centric SND models, our formulation adopts a strategic shipper-oriented viewpoint that removes fixed fleet resource constraints while introducing complex requirements for operational stability and long-haul synchronization. 
    To our knowledge, this is the first formulation in the literature that simultaneously incorporates all these elements at an industrial scale.

    \item[\textbf{(ii) A scalable metaheuristic framework tailored to this setting.}] Methodologically, we develop a novel Iterated Local Search (ILS) metaheuristic specifically designed to navigate the STPP's combinatorial complexity. 
    The framework integrates: (a) specialized large-neighborhood operators adapted to bin-packing-based consolidation, (b) perturbation mechanisms guided by tractable MILP relaxations, and (c) bundle-based decompositions that yield scalable lower bounds. 
    This approach enables the discovery of high-quality solutions for instances involving 700,000 commodities and 1{,}200{,}000 arcs-a scale where standard SND heuristics and exact decomposition methods, such as Column Generation, reach their practical limits.

    \item[\textbf{(iii) A computational study at industrial scale with realized impact.}] Empirically, we demonstrate that the proposed ILS achieves substantial performance improvements over both Renault's legacy planning solutions and state-of-the-art benchmark heuristics. 
    On the largest ``World'' instances, the algorithm delivers solutions revealing a 23.2\% cost-reduction potential relative to previous operational plans. 
    Importantly, the solution has been deployed in production at Renault, where it currently supports weekly planning decisions.
    The framework is primarily used to identify high-value consolidation corridors and support carrier procurement decisions. 
    Based on Renault's internal business evaluations, the planning recommendations produced by the system correspond to an estimated realizable cost savings at approximately €20 million per year. 
    To support reproducibility, we provide an anonymized version of the industrial data used in this study.

    \item[\textbf{(iv) Managerial insights for global inbound logistics.}] Our analysis provides actionable guidance for large-scale industrial shippers: 
    (a) \textit{Consolidation Fidelity:} We prove that explicit 1D bin-packing is a ``step forward'' for realistic planning; relying on standard ``giant container'' relaxations leads to cost distortions of approximately 20\%. 
    (b) \textit{Price of Regularity:} We quantify the trade-off between operational stability and cost, showing that highly regular plans align with practitioner preferences while sacrificing only marginal potential savings. 
    (c) \textit{Strategic Insourcing:} We identify a volume-dependent threshold for outsourcing, demonstrating that for high-volume global networks, in-house planning consistently outperforms 3PL models.
\end{description}

Overall, our work extends the scope of the service network design literature by introducing a shipper-side, bin-packing-explicit formulation and providing the first computational evidence that such problems can be solved effectively at a global industrial scale.

\section{Case study \& Problem Setting} \label{sec:CaseStudy}

This section presents the empirical and operational context for the transportation planning problem studied in this paper. We focus on Renault’s inbound supply network, a global, high-volume, multi-tier logistics system. 
Beyond its industrial relevance, this case study embodies the scientific challenges addressed in this work and summarized in Table~\ref{tab:state_of_the_art}. 
Its scale, heterogeneity, and operational constraints create a large-scale routing and consolidation problem spanning multiple spatial and temporal scales. 
Unlike carrier-centric SND formulations, the shipper-oriented perspective focuses on long-term transportation policies over larger, sparser networks with more commodities, flexible delivery times, detailed consolidation decisions, and stronger commodity-grouping and regularity requirements.
The remainder of this section details the network's topology, flow structure, and planning requirements, beginning with its physical architecture and consolidation logic, followed by key flow characteristics and logistical challenges. 
A detailed and extended version of this case study is available in Appendix~\ref{appendix:case_study}. 
Finally, we present an informal problem statement that synthesizes these elements into a shipper-side transportation planning problem, capturing the scale, constraints, and cost trade-offs specific to this context.

\subsection{Supply Network}

\begin{figure}[!b]
    \begin{subfigure}
        {0.33\textwidth}
        \centering
        \includegraphics[width=0.95\linewidth]{
            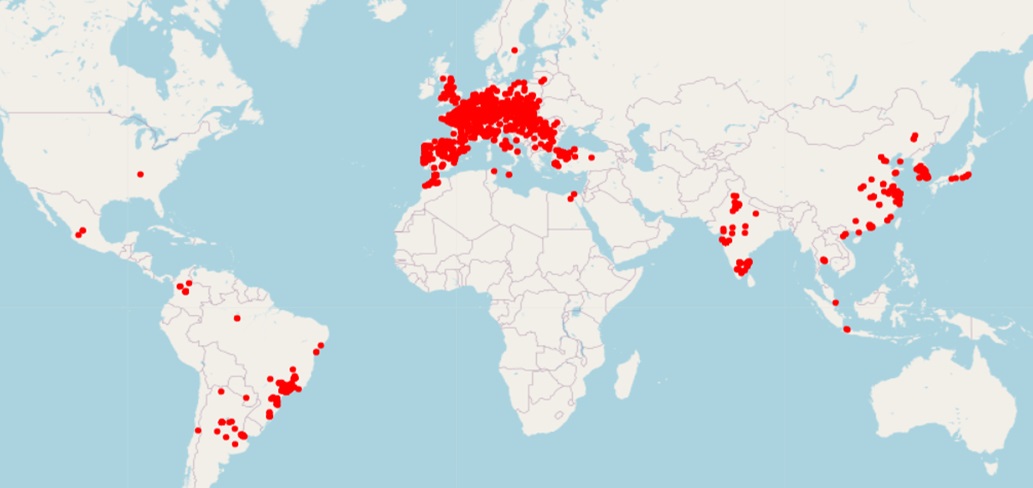
        }

        \smallskip
        \includegraphics[width=0.95\linewidth, height=3cm]{
            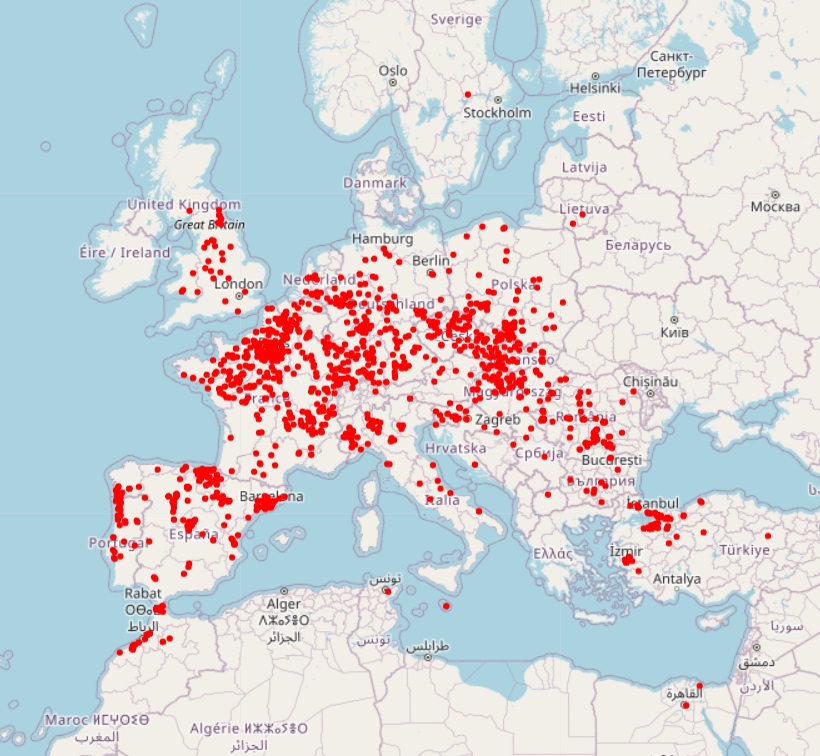
        }
        \caption{\small Suppliers}
    \end{subfigure}
    \begin{subfigure}
        {0.33\textwidth}
        \centering
        \includegraphics[width=0.95\linewidth]{
            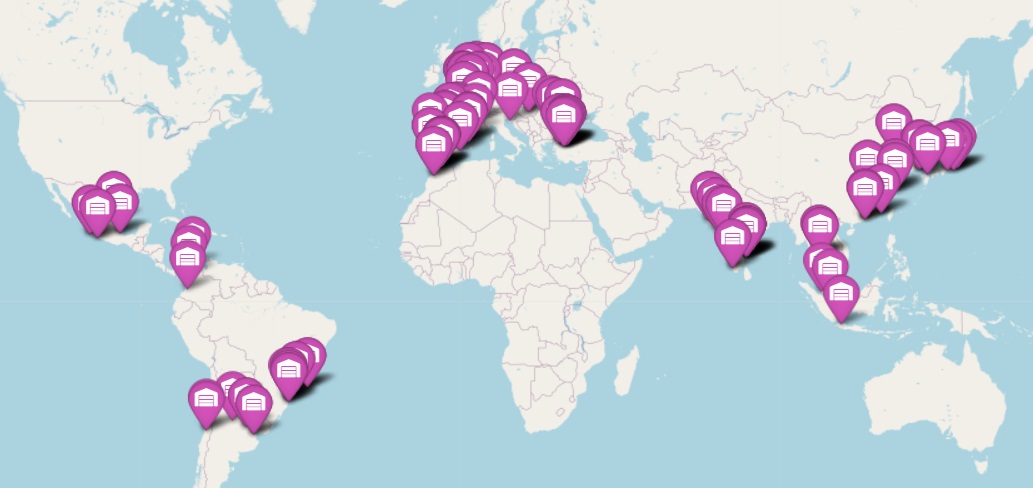
        }

        \smallskip
        \includegraphics[width=0.95\linewidth, height=3cm]{
            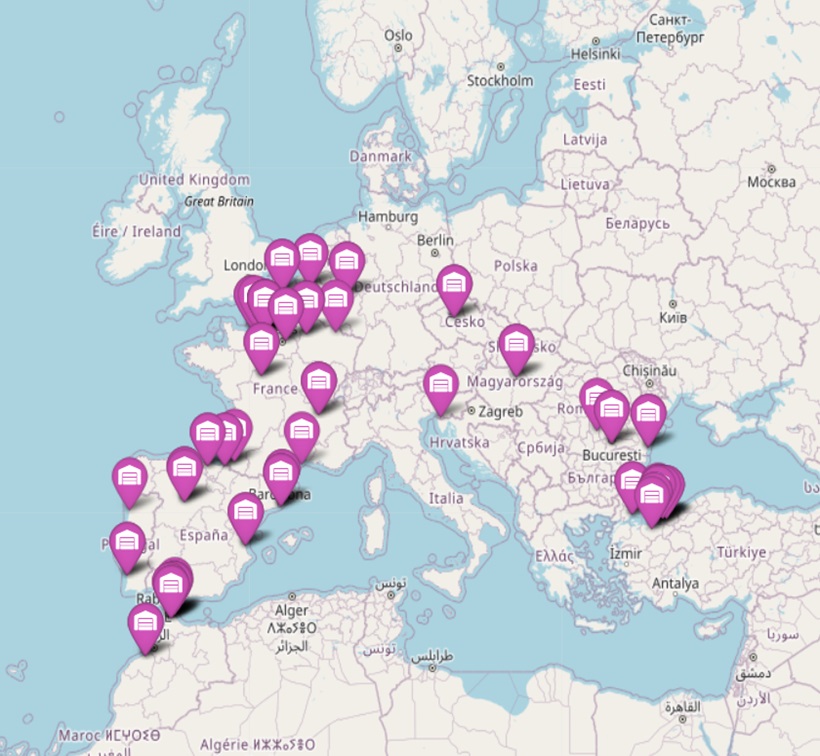
        }
        \caption{\small Platforms}
    \end{subfigure}
    \begin{subfigure}
        {0.33\textwidth}
        \centering
        \includegraphics[width=0.95\linewidth]{
            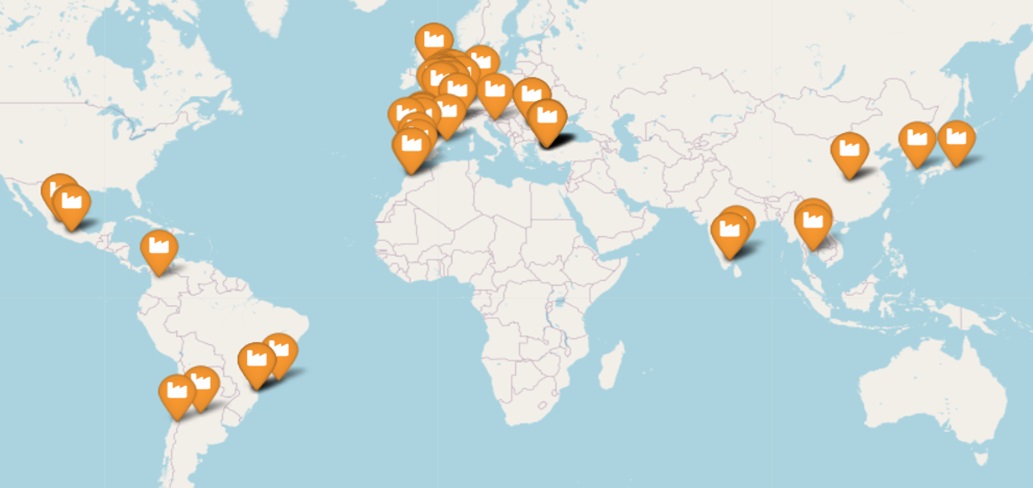
        }

        \smallskip
        \includegraphics[width=0.95\linewidth, height=3cm]{
            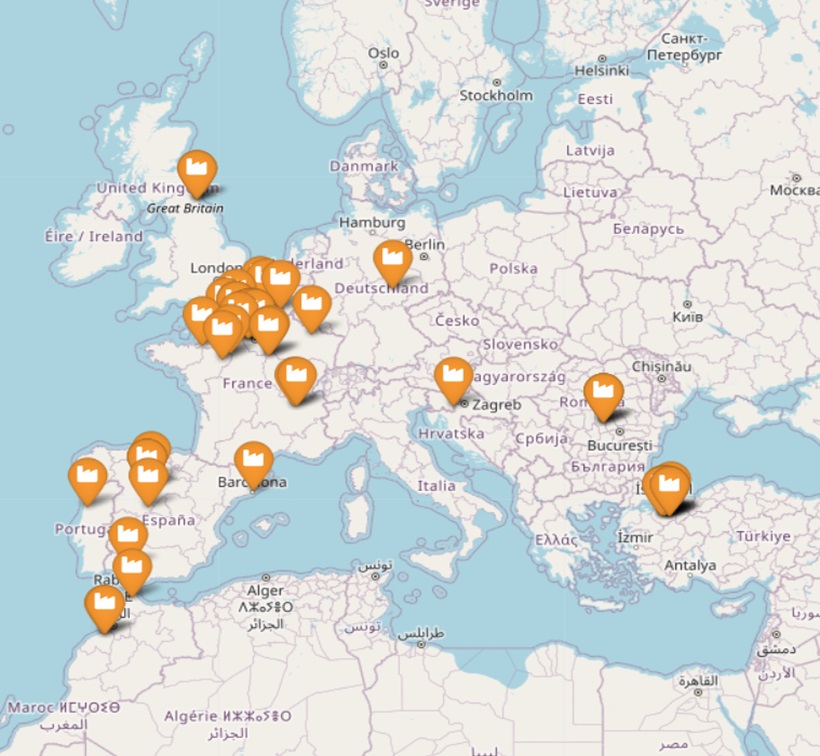
        }
        \caption{\small Industrial Sites}
    \end{subfigure}
    \caption{Worldwide Renault's sites}
    \label{fig:sites-world-1}
\end{figure}

\begin{figure}[!t]
    \begin{subfigure}
        {0.33\textwidth}
        \centering
        \includegraphics[width=0.95\linewidth]{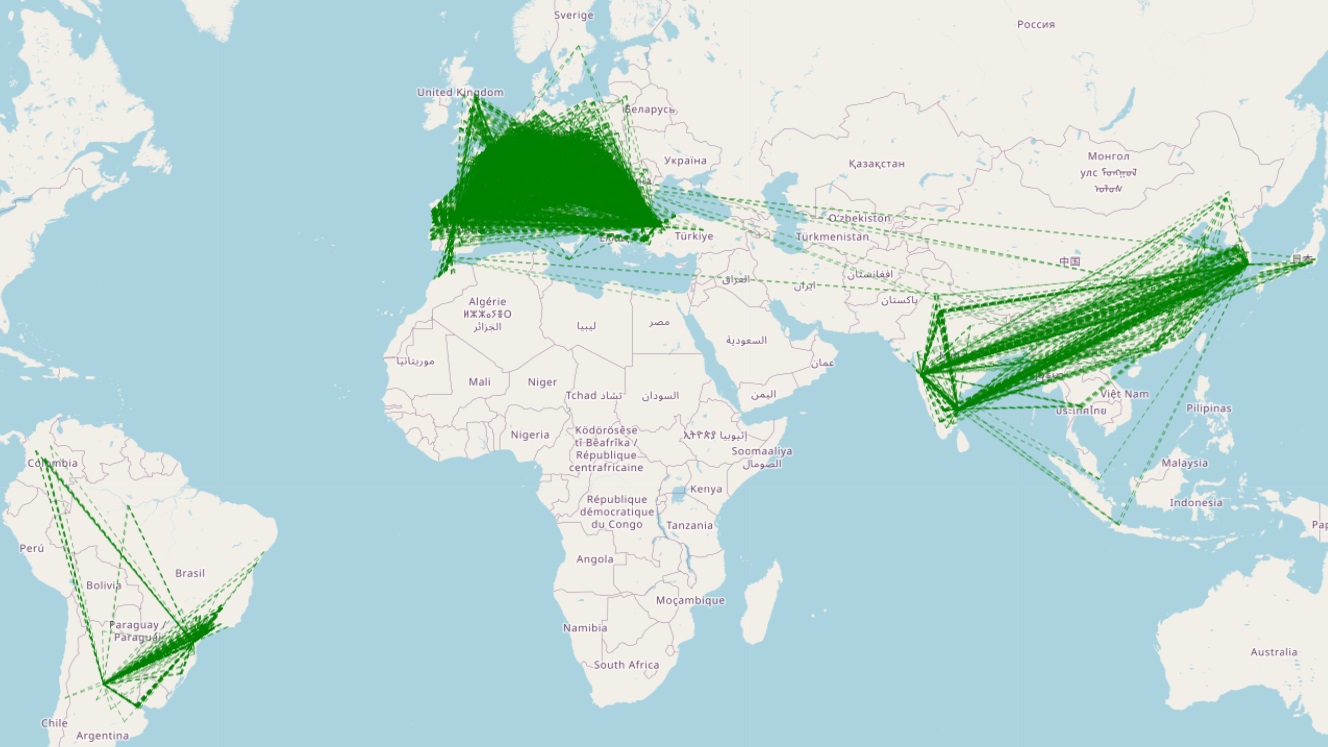}

        \smallskip
        \includegraphics[width=0.95\linewidth, height=3cm]{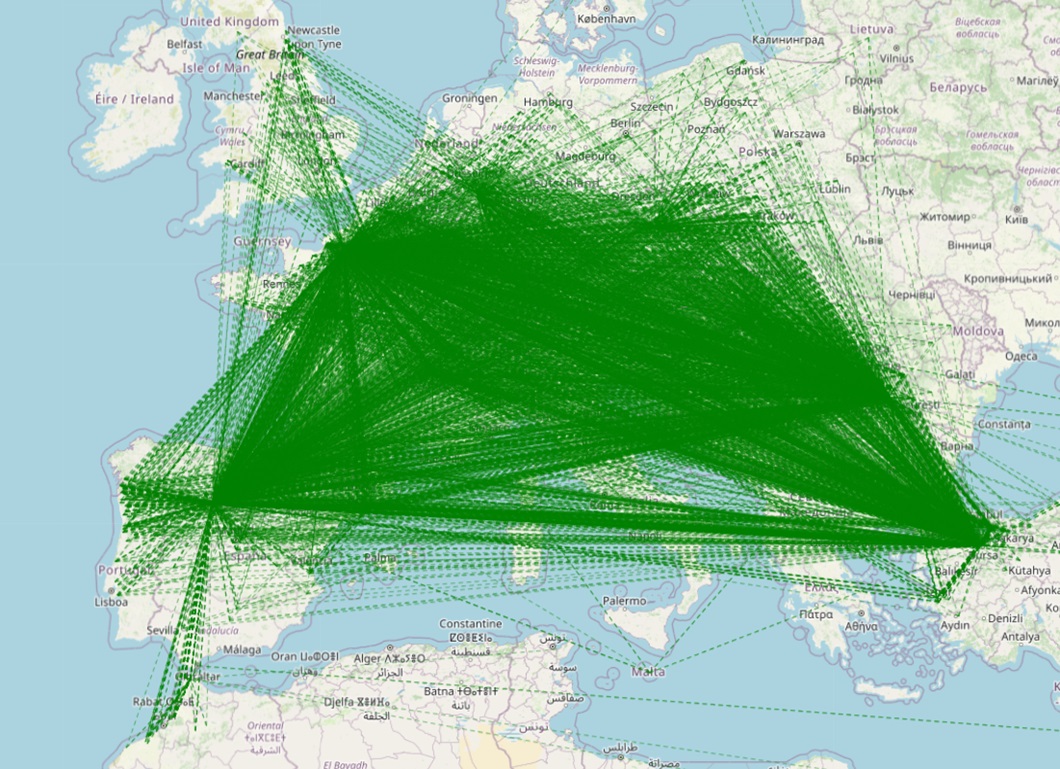}
        \caption{\small Supplier $\to$ Platform}
    \end{subfigure}
    \begin{subfigure}
        {0.33\textwidth}
        \centering
        \includegraphics[width=0.95\linewidth]{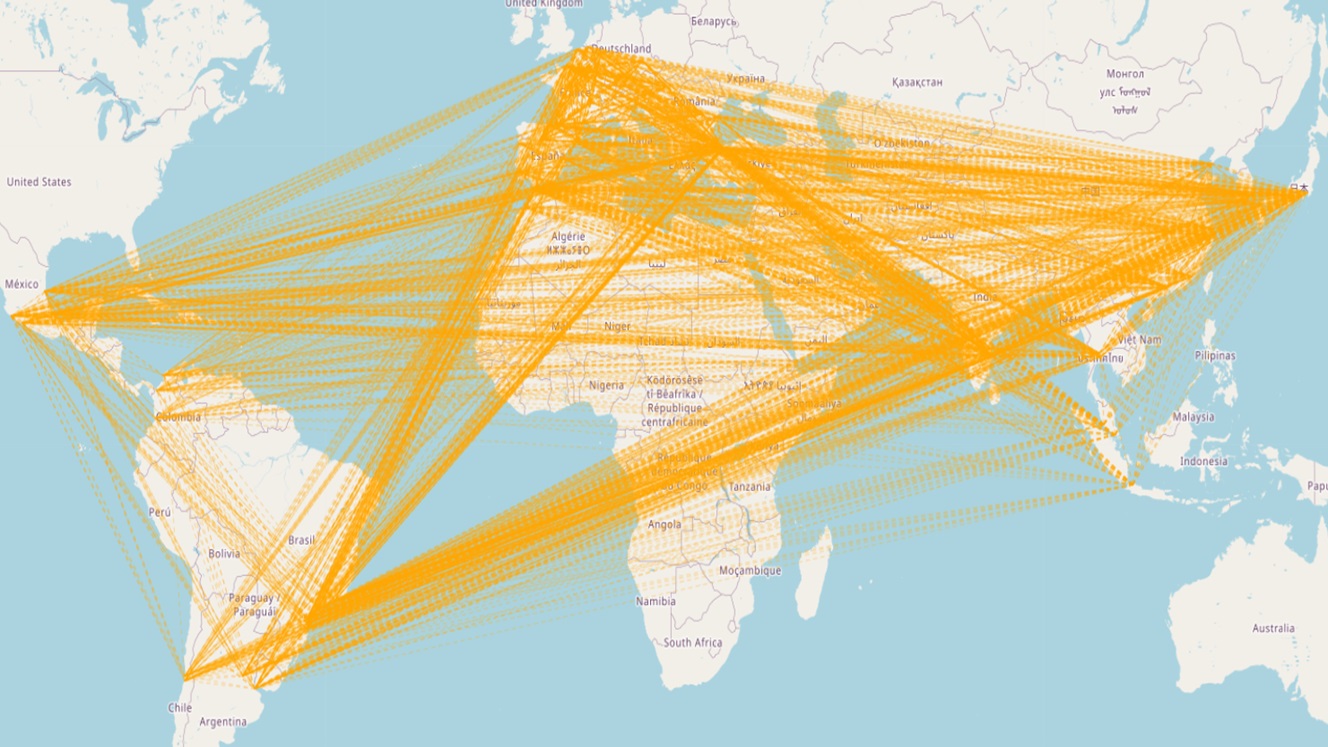}

        \smallskip
        \includegraphics[width=0.95\linewidth, height=3cm]{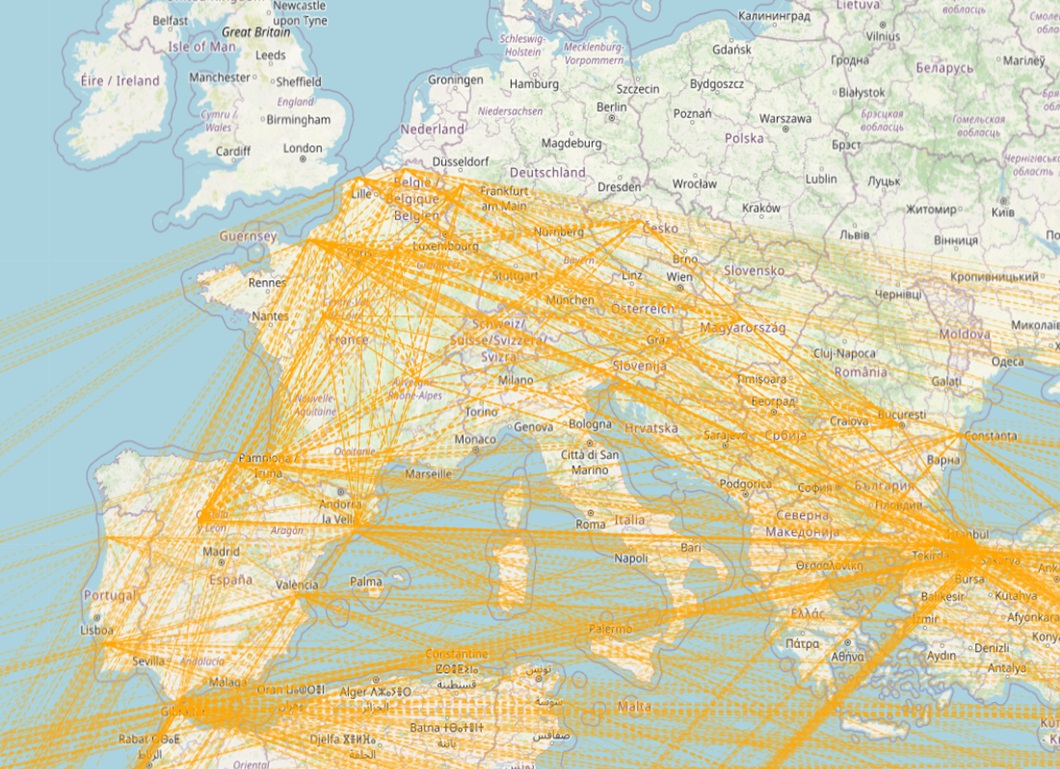}
        \caption{\small Platform $\to$ Platform}
    \end{subfigure}
    \begin{subfigure}
        {0.33\textwidth}
        \centering
        \includegraphics[width=0.95\linewidth]{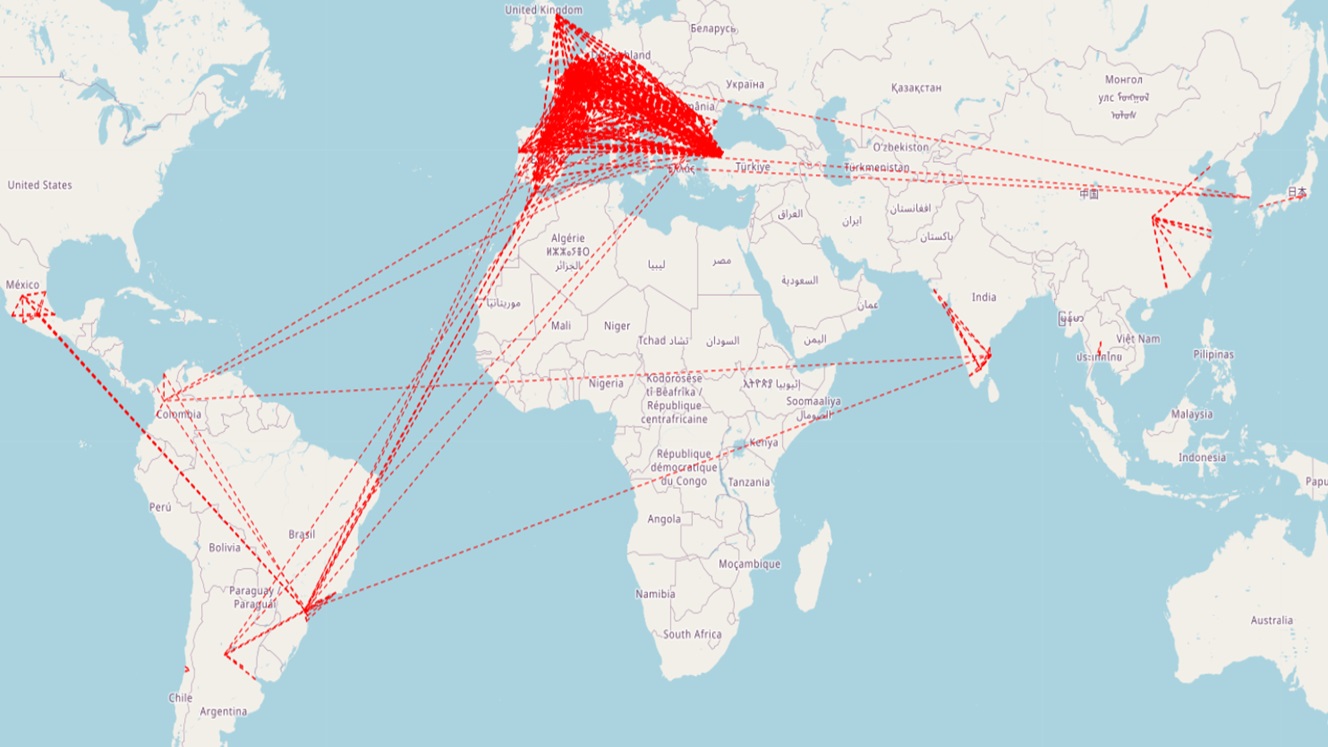}

        \smallskip
        \includegraphics[width=0.95\linewidth, height=3cm]{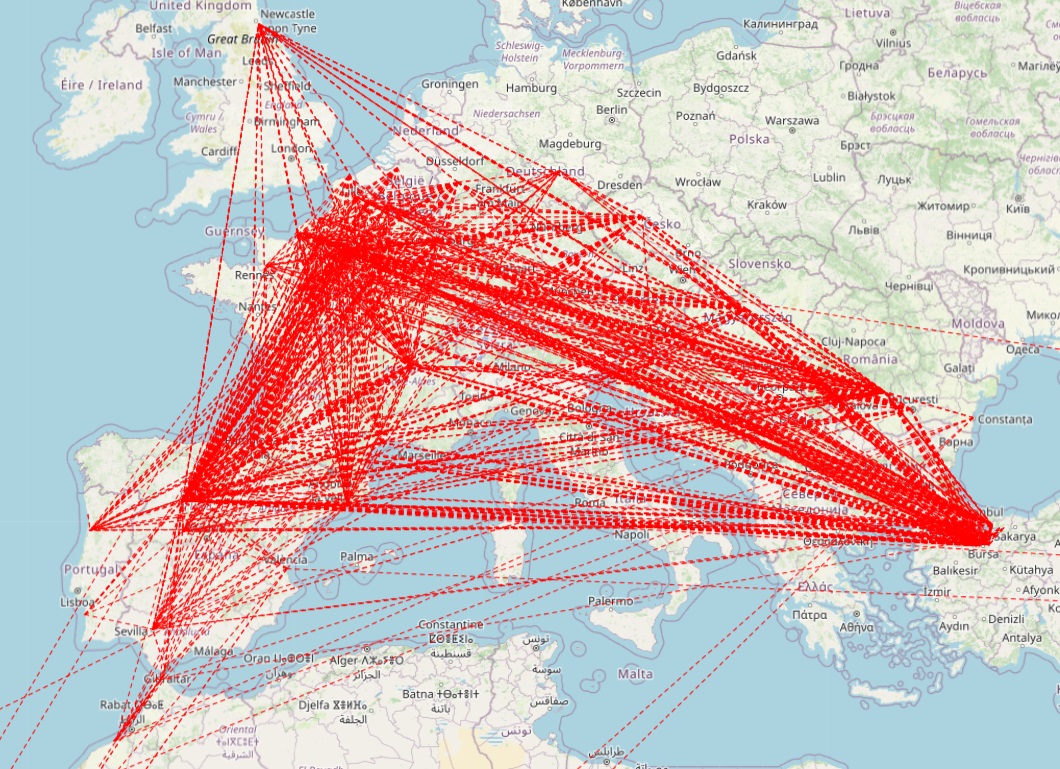}
        \caption{\small Platform $\to$ Sites}
    \end{subfigure}
    \caption{Worldwide Renault's legs}
    \label{fig:legs-world-1}
\end{figure}

Figure~\ref{fig:sites-world-1} illustrates the structure of Renault's supply network.
It includes all suppliers and industrial sites involved in the transport of car parts, as well as all logistics platforms used to consolidate flows.     
We refer to transportation between two nodes as a leg.
Figure~\ref{fig:legs-world-1} visualizes all such legs within the network. 
These maps highlight the global scale and complexity of Renault’s supply network, which comprises more than 3,000 sites, 40 industrial sites, and 100 logistics platforms.
Most sites are located in Europe, and suppliers account for the vast majority of nodes. 
This scale, combined with heterogeneous lead times—ranging from a few days for regional flows to several weeks for intercontinental shipments—creates significant planning complexity, amplified by exposure to disruptions and the need to coordinate consolidation across multiple tiers.
The network follows a hub‑and‑spoke structure: suppliers form a sparsely connected periphery, while platforms act as high‑degree hubs enabling consolidation and reducing the number of direct supplier–plant connections. 
Such configurations support economies of scale while also improve scalability and resilience, as they allow new suppliers to be integrated through local connections to existing hubs and enable rerouting in response to disruptions.
However, fully realizing these benefits requires sophisticated planning tools to manage consolidation, timing, and routing decisions effectively.
Network patterns further reflect this structure. 
Most legs connect suppliers to industrial sites and have short durations; platform‑to‑platform legs, although fewer, are long‑haul, strategically planned, and account for a disproportionate share of transportation cost and emissions. 
Consolidation at platforms is therefore essential for cost efficiency and sustainability, especially for intercontinental flows. 
In this work, each leg is served by a single mode (trucks inland, vessels overseas), reflecting data availability constraints, though the model is fully compatible with a multimodal settings-as discussed in Appendix~\ref{appendix:model_extensions}.
Overall, the network combines a high‑density consolidation core with a vast, low‑connectivity supplier periphery, a topology that drives the methodological requirements motivating the modeling framework developed in the following sections.

\subsection{Commodities and Flows}

\begin{figure}[!b]
    \begin{subfigure}{0.45\textwidth}
        \centering 
        \includegraphics[width=0.95\textwidth, trim= 0pt 20pt 0pt 0pt, clip]{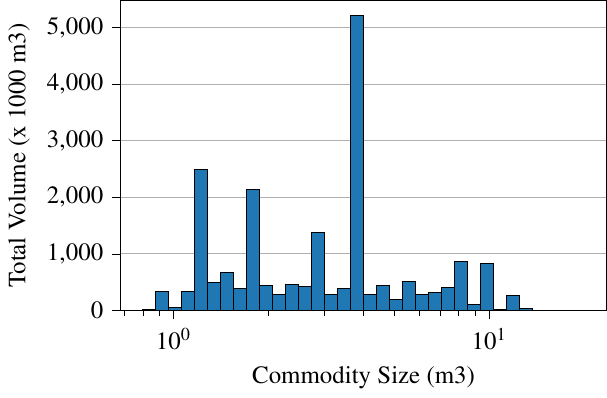}   
        \caption{Total Volume per Size (m$^3$)}
    \end{subfigure}
    \begin{subfigure}{0.45\textwidth}
        \centering 
        \includegraphics[width=0.95\textwidth, trim= 0pt 20pt 0pt 0pt, clip]{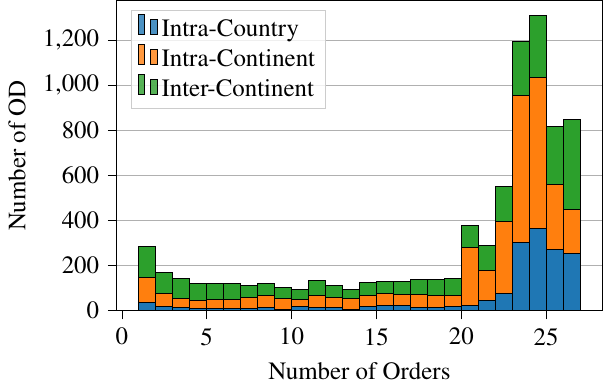}   
        \caption{ Regularity on 26 time steps}
    \end{subfigure}
    \caption{Size and Time Regularity distributions}
    \label{fig:size-order-distribution}
\end{figure}

\begin{figure}[!t]
    \begin{subfigure}
        {0.33\textwidth}
        \centering
        \includegraphics[width=0.95\linewidth]{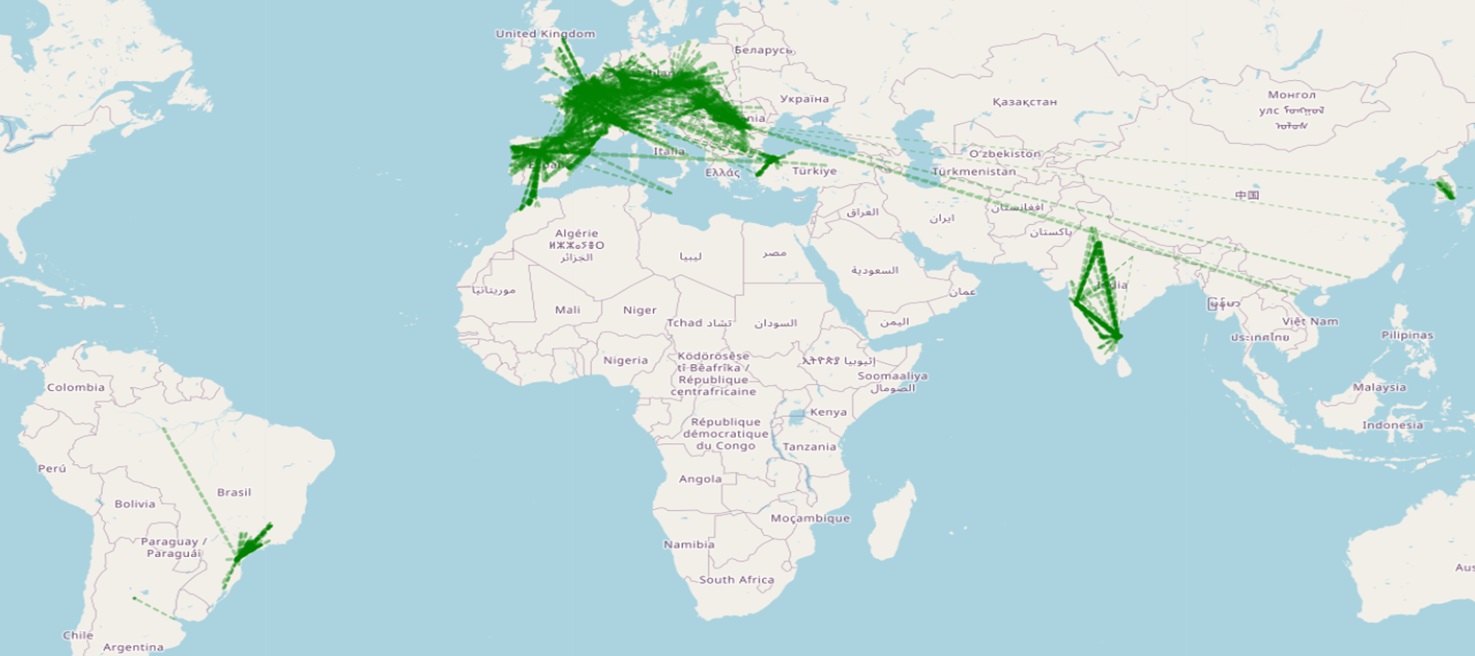}

        \smallskip
        \includegraphics[width=0.95\linewidth, height=3cm]{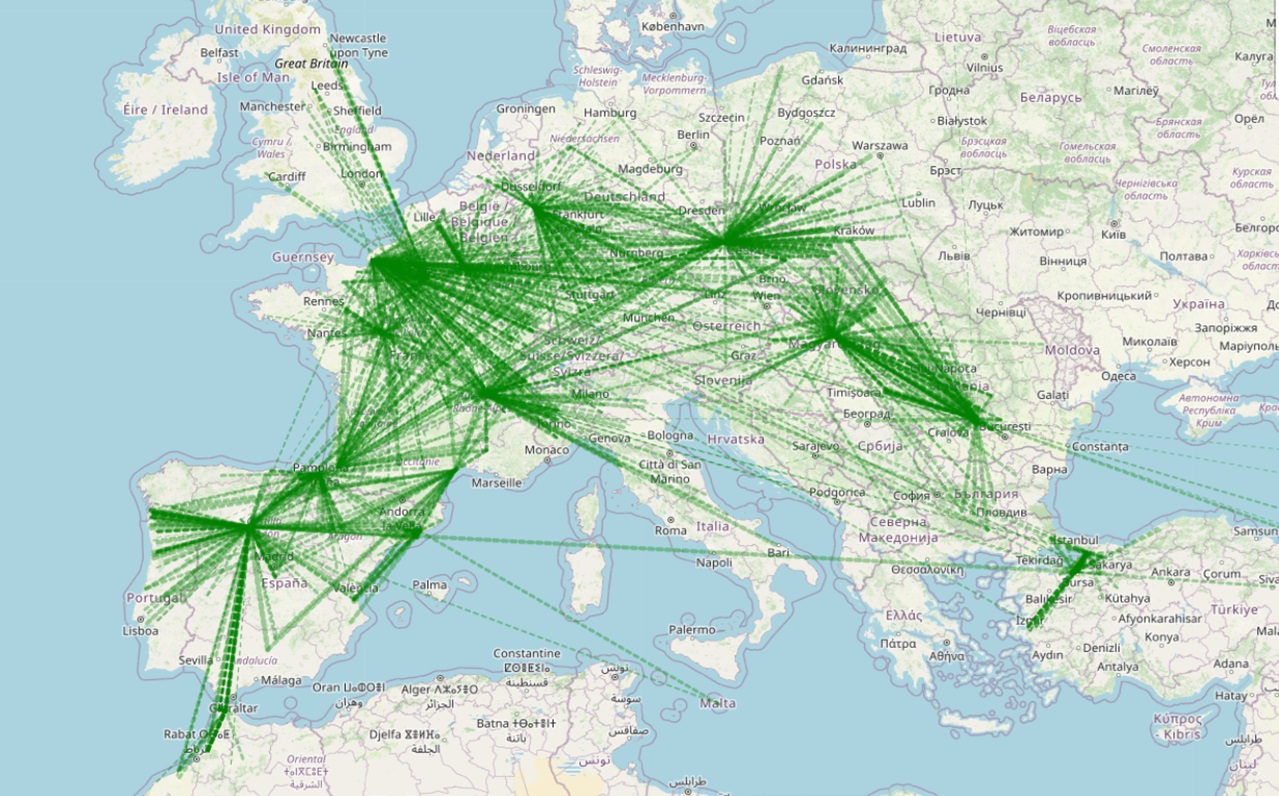}
        \caption{\small Supplier $\to$ Platform}
    \end{subfigure}
    \begin{subfigure}
        {0.33\textwidth}
        \centering
        \includegraphics[width=0.95\linewidth]{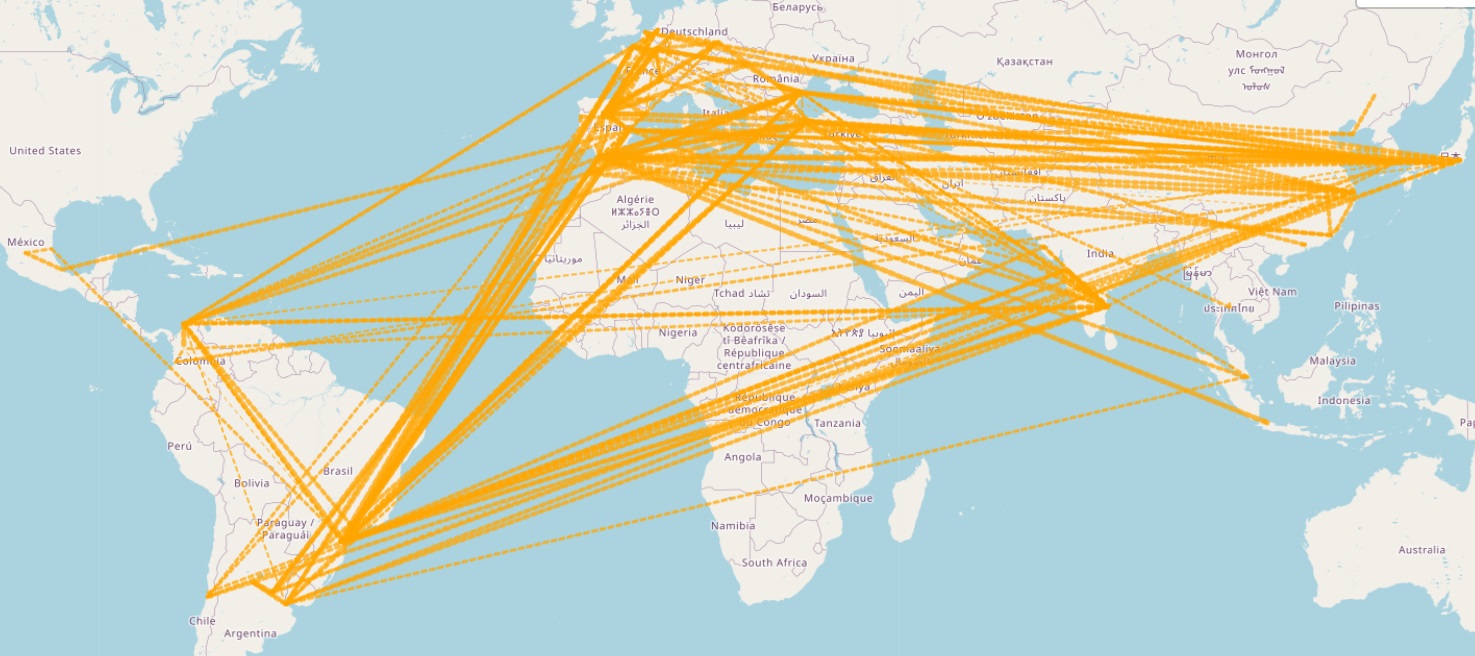}

        \smallskip
        \includegraphics[width=0.95\linewidth, height=3cm]{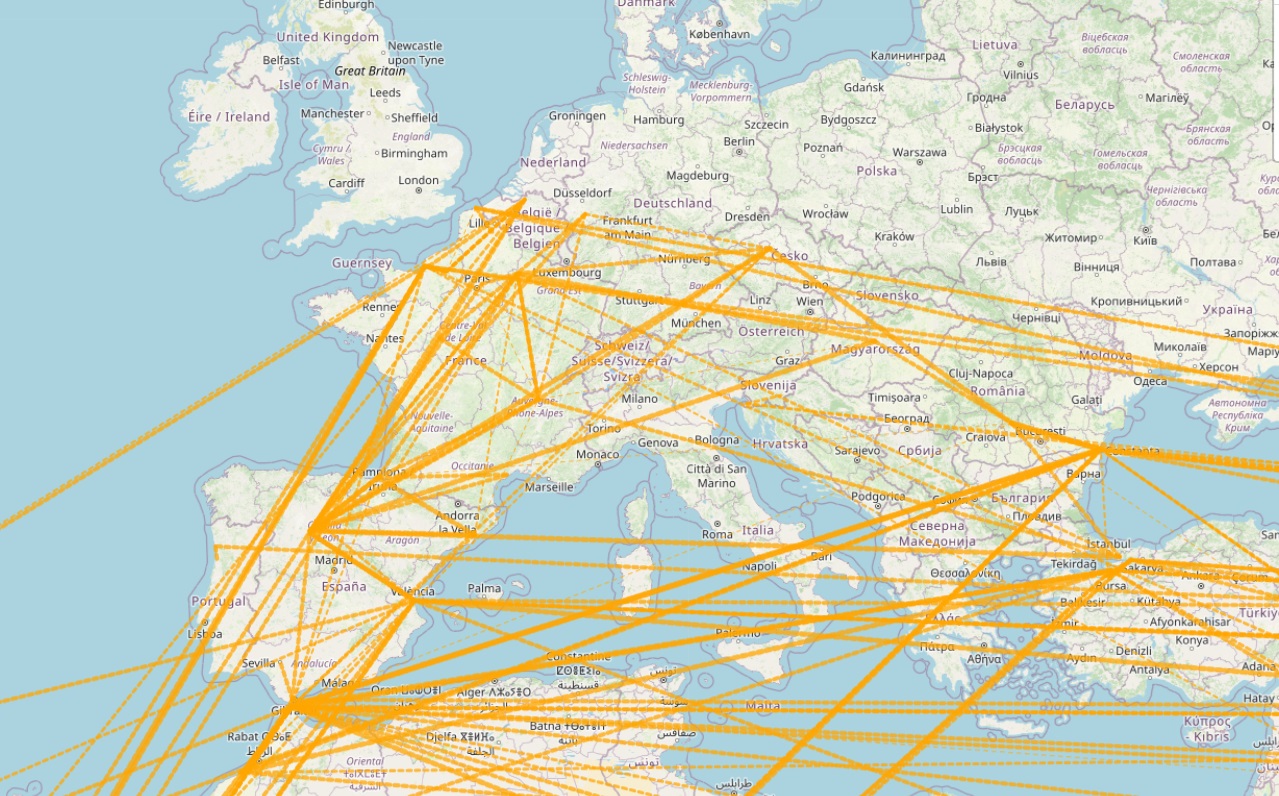}
        \caption{\small Platform $\to$ Platform}
    \end{subfigure}
    \begin{subfigure}
        {0.33\textwidth}
        \centering
        \includegraphics[width=0.95\linewidth]{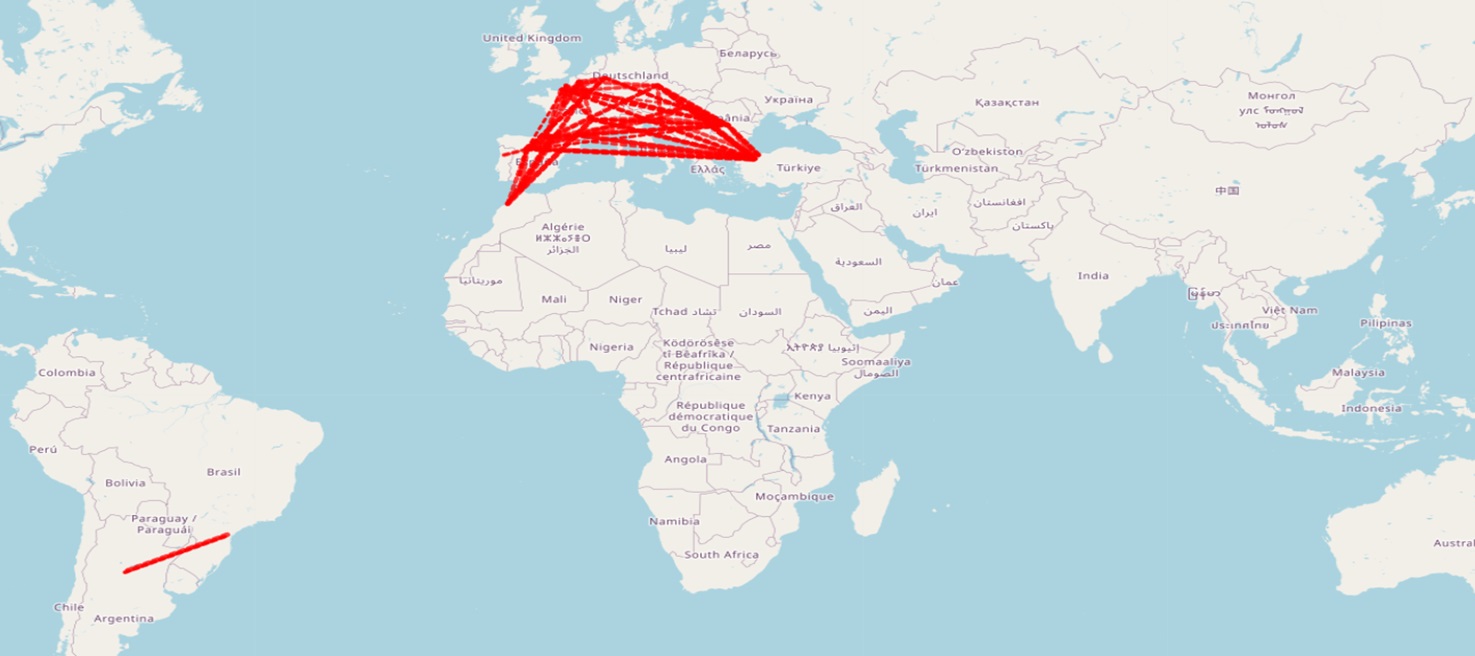}

        \smallskip
        \includegraphics[width=0.95\linewidth, height=3cm]{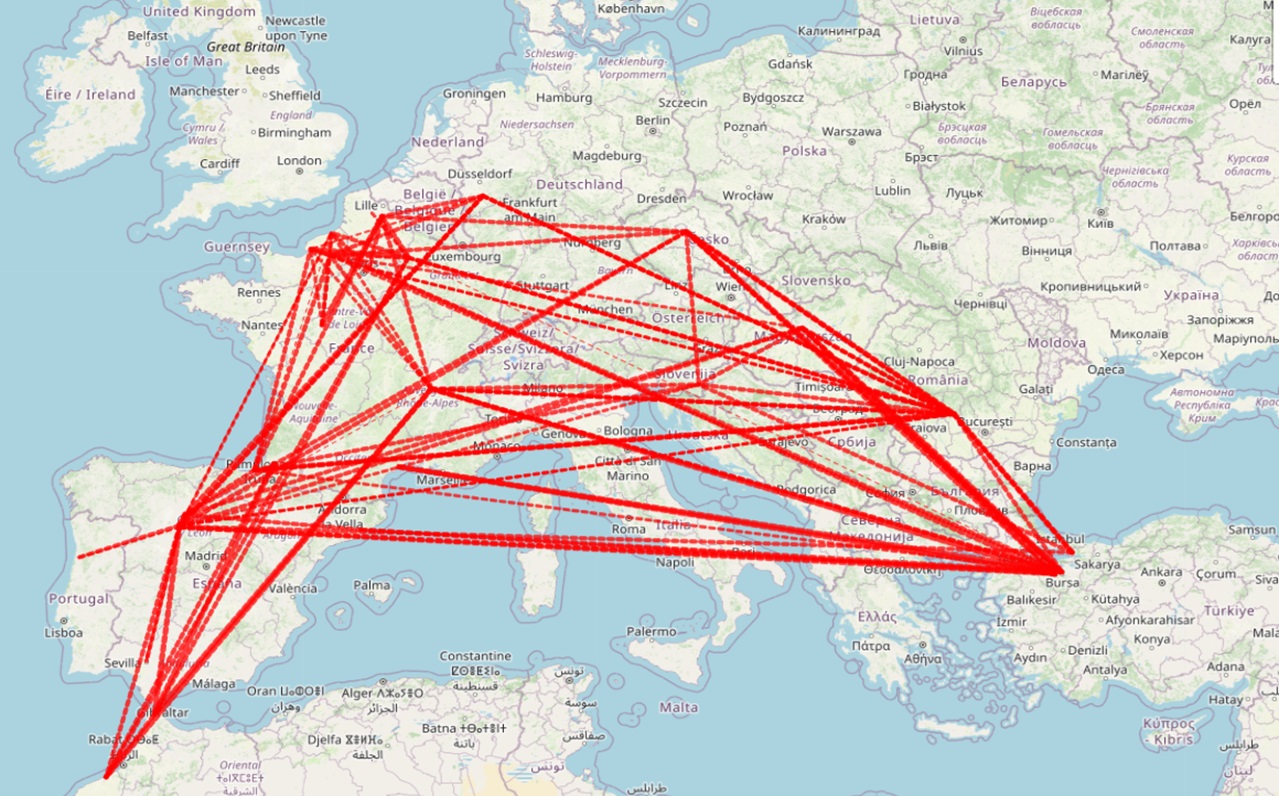}
        \caption{\small Platform $\to$ Sites}
    \end{subfigure}
    \caption{Legs used in current flows}
    \label{fig:flows-world-1}
\end{figure}

Renault’s inbound supply chain supports a highly complex and large-scale operation.
Over a six-month horizon, the network handles more than 40,000 distinct car parts, divided into 700,000 commodities to deliver, corresponding to approximately 9,000,000 packages and 20,000,000 m$^{3}$ of volume.
Figure~\ref{fig:size-order-distribution} illustrates the distribution of this total volume across shipment sizes and demand regularity. 
Despite this magnitude, inbound flows remain relatively stable over time, enabling structured mid‑term planning. 
Still, the system must remain adaptable to changes in supplier configurations or demand.
In this context, the sheer number of parts and the high granularity of the commodity mix creates substantial planning complexity, as parts differ widely in size, frequency, and origin–destination patterns. 
Most shipments fall in the 1–4 m$^3$ range and cannot fill a service alone, making efficient consolidation essential. 
This heterogeneity increases the dimensionality of planning problems and necessitates sophisticated aggregation, routing, and consolidation strategies to efficiently organize transportation, especially for multi‑leg paths.
Across the network, nearly 500,000 transport services—trucks or shipping containers—are used every six months, with Figure~\ref{fig:flows-world-1} illustrating their spatial distribution.

Yet only a small fraction of potential legs is activated. 
Flows are funneled through a limited set of high‑volume corridors where consolidation is most effective.
Finally, while most physical volume moves within Europe, intercontinental legs—though fewer—account for a disproportionate share of cost and emissions.
Overall, Renault’s inbound flow structure combines high commodity granularity, relative temporal stability and selective leg activation. 
Effective planning must therefore leverage consolidation opportunities, operate across multiple temporal and spatial scales, and dynamically determine which legs and routes to activate. 

\subsection{Planning problem}

As stated in the introduction, Renault’s transportation planning process works at the strategic level and aims to generate a provisional six-month plan for delivering car parts from suppliers to assembly industrial sites across its global logistics network.
The case study illustrates several of the structural characteristics that define the scientific challenges investigated in this work. 
It introduces a large-scale transportation planning problem combining routing, scheduling, and packing considerations.
The resulting task belongs to the class of Shipper Transportation Planning Problems, which we may describe informally as follows:
\begin{equation}\label{eq:PbInformal}
   \begin{aligned}
       \min \, & \sum_{a} \left (
            \begin{array}{c}
              \text{ transport cost (€/service) } \\
              \text{ platform cost (€/m$^3$) } \\
              \text{ capital cost (€/km) }
            \end{array}
            \right )_a \\
        \mathrm{s.t.} \, & \left |
\begin{array}{l}
  \text{ Routes with regularity requirements} \\
  \text{ Flexible delivery time} \\
  \text{ Bin-packing consolidation inside services}
\end{array}
\right .
   \end{aligned}
\end{equation}

At its core, the planning problem consists of routing a large number of commodities through a multi‑tier logistics network while satisfying operational constraints. 
The objective function reflects Renault’s multi‑dimensional cost structure, combining per‑service transport costs, volume‑based handling costs at platforms, and distance‑based capital‑in‑transit costs. 
These components capture the trade‑off between direct shipments and consolidation through intermediate platforms.
The constraints mirror the operational realities of the network. 
Flows must respect feasible routes, available modes, and delivery timing, with flexible delivery windows allowing for typical transport variability. 
Regularity constraints promote stable routing patterns across time, reflecting planners’ preferences.
Explicit bin‑packing constraints model consolidation decisions at the service level—an important requirement for Renault because the diversity of part volumes makes simple capacity approximations inadequate for high‑quality planning.

The role of the model within Renault's planning process is summarized in Table~\ref{tab:decision_layers}. 
It is primarily intended to support medium-term managerial decisions, including carrier procurement and the identification of major consolidation corridors. 
\begin{table}[t]
    \centering
    \caption{Decision layers supported by the proposed transportation planning model.}
    \label{tab:decision_layers}
    \renewcommand{\arraystretch}{1.3} 
    \resizebox{\textwidth}{!}{%
        \begin{tabular}{p{0.4\textwidth} p{0.4\textwidth} p{0.4\textwidth}} 
            \toprule 
            \textbf{Supported Decisions} & \textbf{Partially Supported Decisions} & \textbf{Unsupported Decisions} \\ 
            \midrule
            Carrier contracts & Network evolution & Specific service bookings \\ 
            Capacity negotiations and procurement & supplier-to-platform policy evaluation & Detailed vehicle load plans \\ 
            Strategic consolidation-corridor selection & Recurring routing structures &  Day-to-day transportation execution \\ 
            \bottomrule 
        \end{tabular}
    }
\end{table}
From this perspective, routing and bin-packing decisions serve as planning proxies that estimate how freight is likely to consolidate across the network and therefore determine the true economic attractiveness of candidate transportation corridors. 
Explicitly modeling consolidation remains essential because strategic procurement decisions depend directly on realistic estimates of vehicle utilization and transportation costs as supported by the managerial analyses in Section~\ref{subsec:ext_analysis}.
The decision model also supports robust managerial recommendations. 
The optimization is repeatedly updated as demand forecasts evolve, while regularity constraints favor recurring transportation structures over isolated short-term opportunities.
The formulation therefore integrates the structural properties of automotive inbound supply chains’s logistics system with its key operational constraints. 
The remainder of the paper develops a scalable algorithmic framework capable of addressing this complexity.

\section{Encoding structure and regularity through time expanded graphs} \label{sec:PbStatement}

The transportation planning problem faced by Renault involves a highly complex decision space due to the scale of operations and the rich structure of operational constraints.
Each planning instance involves millions of part-level shipments over a multi-tier network with temporal dependencies.
In addition to traditional network flow constraints, the problem introduces explicit bin-packing requirements for transport services, coupling flow-assignments non-linearly, which substantially increase the computational complexity.
Furthermore, classical service network design models typically rely on a flat network representation, combined with a predetermined set of admissible paths for each commodity and a time‑expanded graph to map these paths across the planning horizon. 
This modeling framework, however, does not suit our setting: Renault’s inbound system does not provide a predefined path set, generating all feasible paths for each bundle would be computationally prohibitive, and classical expansions do not naturally capture flexible delivery times.
To address these challenges, we develop a graph-based modeling framework that integrates network structure, temporal dynamics, regularity requirements and bin-level consolidation decisions in a scalable and modular way.
Beyond providing a formal problem formulation, the structure of this modeling framework eases the development of a practically efficient algorithm in the subsequent section.

Specifically, we present a formal model composed of four main components: (i) the static supply network and the temporal structure of commodities in Section~\ref{subsec:supply_network_and_commodities}, (ii) two coupled time-expanded graphs to represent delivery timing and regularity in Section~\ref{subsec:time_expanded_graphs}, (iii) the flow, regularity, and bin-packing constraints in Section~\ref{subsec:flow_regularity_bp_contraints}, and (iv) the associated cost structure and resulting mixed-integer programming (MIP) formulation in Section~\ref{subsec:cost_and_full_pb}.
We conclude with a discussion of model scalability and possible extensions.

\subsection{Modeling the Supply Network and Commodities}\label{subsec:supply_network_and_commodities}

\begin{figure}[!hb]
    \centering
    \includegraphics[scale=0.9]{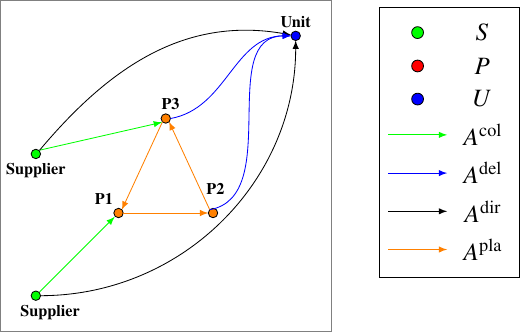}
    \caption{Example of network graph $D$.} 
    \label{fig:flatNetworkExample}
\end{figure}

We model the supply network as a directed graph $D = (V, A)$, where the node set $V = S \cup P \cup U$ includes all physical locations in the network.
The set $S$ contains suppliers, $P$ represents logistics platforms, and $U$ denotes production units.
The arc set $A \subseteq V^2$ captures all possible transportation links between locations and is partitioned into four categories.
\emph{Collection} arcs $A^{\mathrm{col}}$ connect suppliers to platforms and represent the first leg of multi-leg paths.
\emph{Inter-platform arcs} $A^{\mathrm{pla}}$ connect different platforms and enable the transfer of consolidated goods within the intermediate logistics network.
\emph{Delivery arcs} $A^{\mathrm{del}}$ link platforms to production units and represent the final leg of consolidated flows.
Finally, \emph{direct arcs} $A^{\mathrm{dir}}$ connect suppliers directly to production units, bypassing intermediate consolidation.
Figure~\ref{fig:flatNetworkExample} illustrates a representative example of such a network, comprising two suppliers, three platforms, and one production unit.

The inbound supply chain operates on a rolling time horizon $[T] = {1, \ldots, T}$, discretized in our application-specific context into weekly steps.
This rolling representation reflects the assumption that future transportation needs exhibit sufficient regularity.
A \emph{commodity} $m \in M$ represents a non-splittable part sent from supplier $s_m \in S$ to unit $u_m \in U$ for delivery in week $d_m \in [T]$.
Each commodity has a volume $\ell_m \in \mathbb{R}$, a maximum delivery time $\tau_m \in \mathbb{N}$ (in weeks), and a multiplicity $q_m$ representing the quantity to be delivered in week $d_m$.
The specific delivery path used by each commodity is a decision variable.
In this context, we note that we treat each supplier–unit–delivery tuple as a distinct commodity, although the same part may be sourced from multiple suppliers or sent to multiple industrial sites.
This simplification aligns with the scope of this work, which focuses on transportation planning rather than sourcing decisions, as this is standard modeling practice in the transportation planning literature \citep{crainic202550}.

To structure delivery decisions and ensure operational regularity, we introduce two key constraints: grouping and time regularity.
\emph{Grouping} enforces that all commodities sharing the same origin supplier, destination unit, and delivery date must follow an identical path through the network.
To formalize this, we aggregate such commodities into \emph{orders} $o \subseteq M$.
Each order $o$ is characterized by a common delivery date $d_o$, a shared maximum delivery time $\tau_o$, and a total volume $\ell_o = \sum_{m \in o} q_m \ell_m$.
Let $O$ denote the set of all orders.
\emph{Time regularity} imposes consistency in routing over time: commodities with identical supplier–unit pairs must follow the same sequence of nodes in the supply network, regardless of their delivery date.
To capture this, we group all orders with the same supplier and unit but differing delivery dates into \emph{bundles} $b \subseteq O$, which represent recurring flows along regular paths.
Let $B$ denote the set of all bundles.
This regularity structure directly reflects the operational realities observed in the Renault case study: grouping promotes consolidation of freight, while temporal regularity promotes repeatable transportation plans that can be implemented by logistics operators.

Together, these two layers of aggregation form a hierarchical structure: the set of orders $O$ partitions the set of commodities $M$, and the set of bundles $B$ partitions the set of orders $O$.
This structure ensures that all commodities within a bundle follow the same node sequence, thereby enforcing both intra-week consistency (via grouping) and inter-week regularity (via bundling) in delivery planning.

\subsection{Time-Expanded Graph Structures}\label{subsec:time_expanded_graphs}

\begin{figure}[!b]
    \centering
    \includegraphics[width=0.9\textwidth, height=0.65\textheight]{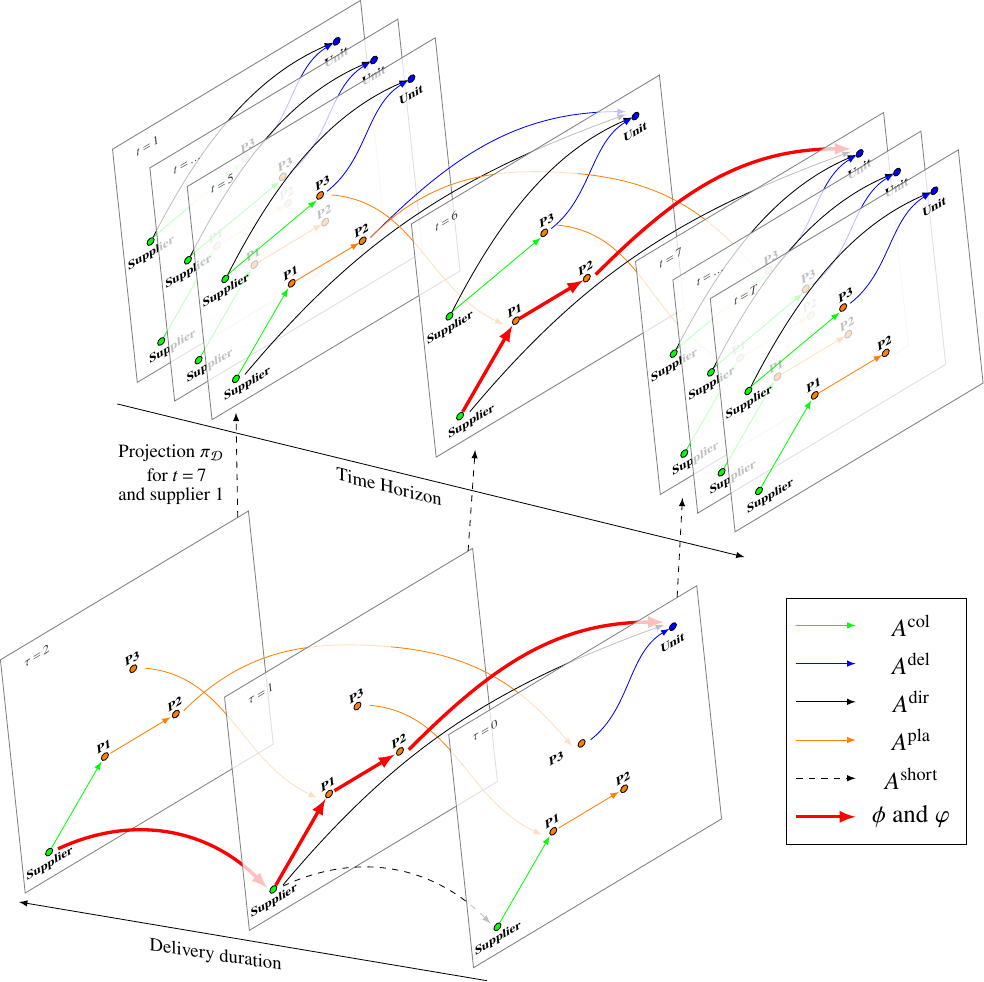}
    \caption{Example of time-space digraph $\mathcal{D}$ (top) and travel-time digraph $\mathscr{D}$ (bottom)} 
    \label{fig:travelTimeExample}
    \vspace{-0.25cm}
    \note{On the travel-time graph $\mathscr{D}$ (bottom), given a bundle $b \in B$, we compute a single path $\phi_b$, thereby enforcing regularity requirements across all associated orders.
    For each order $o \in b$, we translates $\phi_b$ into $\varphi_o = \pi_{\mathcal{D}}(\phi_b, o)$ on $\mathcal{D}$ (top graph) according to its delivery date recovering the actual timed path used by the commodities $m \in o$.}
\end{figure}

To model both the temporal flow of commodities and the regularity of delivery paths, we define two complementary time-expanded graphs.
While time-expanded graphs are a standard modeling tool, the simultaneous enforcement of shipment regularity constraints and flexible delivery times would lead to prohibitively large formulations if handled directly on a single time-space representation. 
To address this challenge, we introduce a second, travel-time-based graph together with a projection mechanism linking both structures. 
This design allows regularity constraints to be represented directly within the graph topology while shifting most path computations to a significantly smaller graph.
Figure~\ref{fig:travelTimeExample} illustrates these two time expansions and their relations, while purposely omitting rolling-horizon forming arcs for clarity.

\paragraph{Time-Space Graph}

To capture the temporal dynamics of the supply network, we first construct the classical time-space graph.
We expand the static graph $D$ into a time-indexed structure that reflects both fluctuating demand and transportation durations over the planning horizon.
This graph provides the reference representation on which transportation flows, delivery dates, and transport durations are explicitly tracked over the planning horizon.
This results in the \emph{time-space graph}, denoted by $\mathcal{D} = (\mathcal{V}, \mathcal{A})$.
We define the set of \emph{timed nodes} as $\mathcal{V} = V \times [T]$, where each node $\nu = (v, t)$ represents location $v \in V$ at time step $t \in [T]$.
For each arc $a = (u, v) \in A$, we let $\tau_a \in \mathbb{N}$ denote the number of time steps required to traverse the arc.
Using this, we define the set of \emph{timed arcs} as
$$\mathcal{A} = \left\{ \big((u, t), (v, t')\big) \;\middle|\; (u, v) \in A \; \text{and} \; t' \equiv t + \tau_{(u,v)} \pmod{T} \right\} \subseteq \mathcal{V}^{2}$$
The modular operator rolls the time horizon, allowing arcs that start near the end of the horizon to wrap around to the beginning of the next cycle.
We use the symbol $a$ for arcs in either $A$ or $\mathcal{A}$, relying on context to distinguish between the two.

\paragraph{Travel-Time Graph}

To model regularity constraints and account for flexible delivery times, we define the \emph{travel-time graph} $\mathscr{D} = (\mathscr{V}, \mathscr{A})$.
This graph constitutes the main conceptual contribution of our modeling approach. 
Rather than computing transportation paths independently for every order on the time-space graph and then enforcing regularity constraints, regularity is embedded directly into the graph structure. 
A single path computed on $\mathscr{D}$ simultaneously defines the routing decisions of all orders belonging to the same bundle, thereby transforming a large set of coordination constraints into a graph projection.
This construction provides a computational advantage: routing decisions are represented through network flows on a smaller graph spanning the delivery-time horizon rather than the full planning horizon, preserving the favorable structure of network-flow formulations that is later exploited within several neighborhood optimization subproblems.
This graph is a partial time-expansion of the network graph $D$, where time steps represent the remaining time until delivery rather than absolute positions on the planning horizon.
Unlike the time-space graph, which expands over the full horizon $[T]$, the travel-time graph expands over the partial horizon $[\mathcal{T}]$, where $\displaystyle \mathcal{T} = \max_{o \in O} \tau_o$ denotes the maximum delivery time allowed for any order.
We expand suppliers and platforms across all time steps in $[\mathcal{T}]$, while production units only appear at the final time step $\mathcal{T}$.
Accordingly, we define the node set as $\mathscr{V} = \mathscr{S} \cup \mathscr{P} \cup \mathscr{U}$, where $\mathscr{S} = S \times [\mathcal{T}]$, $\mathscr{P} = P \times [\mathcal{T}]$, and $\mathscr{U} = U \times \{\mathcal{T}\}$.
We construct the arc set $\mathscr{A}$ analogously to that of the time-space graph, with one key extension: we introduce \emph{shortcut arcs} to allow delivery time flexibility.
These arcs are defined as $\alpha \in \mathscr{A}^{\mathrm{short}} = \{ ((s, t), (s, t - 1)) \mid (s, t) \in \mathscr{S} \}$ and represent the option to delay the dispatch of a commodity by one time unit at its origin.

\paragraph{Interactions}

The travel-time graph $\mathscr{D}$ interacts with the time-space graph $\mathcal{D}$ through an order-dependent projection operator $\pi_{\mathcal{D}}$, which maps arcs in $\mathscr{D}$ to corresponding arcs in $\mathcal{D}$.
Let \( \rho_T(k)=1+\big((k-1)\bmod T\big) \) denote the canonical mapping from an absolute time index $k \in \mathbb{Z}$ to its representative in $[T]$.
Consider an order $o \in O$ with delivery date $d_o \in [T]$ and a non-shortcut arc \( \alpha=\big((u,r),(v,r')\big) \), where $r$ and $r'$ denote the numbers of time steps remaining before delivery. 
The corresponding departure and arrival dates on the unrolled time axis are respectively $d_o-r$ and $d_o-r'$. 
The projection of $\alpha$ onto the cyclic time-space graph is therefore defined by \( \pi_{\mathcal{D}}(\alpha,o) = \left( \big(u,\rho_T(d_o-r)\big), \big(v,\rho_T(d_o-r')\big) \right) \).
Shortcut arcs do not represent movements on the physical transportation network, it changes the relative-time index of the first physical arc selected for the order and thus postpones its projected dispatch date by one time step.
This projection mechanism enables a consistent translation between the two graph representations.

\subsection{Flow, Regularity, and Bin-packing constraints}\label{subsec:flow_regularity_bp_contraints}

\paragraph{Flows and Regularity Constraints}
We denote by $f_{a}^{m} \in \mathbb{N}$ the quantity of commodity $m \in M$ flowing on arc $a \in \mathcal{A}$.
We denote by $x_{\alpha}^{b} \in \{0, 1\}$ the binary variable indicating whether bundle $b \in B$ uses arc $\alpha \in \mathscr{A}$.
The flow and regularity constraints define an elementary path on the travel-time graph, as captured by equations~\eqref{constr:bundlePath} and~\eqref{constr:pathElementarity}, and project the corresponding quantities onto the time-space graph, as shown in equation~\eqref{constr:projTTonTS}.
For $b \in B$ and $\nu \in \mathscr{V}$, $e_{\nu}^{b} = 1$ for $\nu = (s_{b}, \mathcal{T}- \tau_{b})$, $-1$ for $\nu = (u_{b}, \mathcal{T})$ and $0$ otherwise.
    \begin{equation}
    \sum_{\alpha \in \delta^+(\nu)}x_{\alpha}^{b}- \sum_{\alpha \in \delta^-(\nu)}
    x_{\alpha}^{b}= e_{\nu}^{b}\quad \forall b \in B \text{, }\nu \in \mathscr{V}
    \label{constr:bundlePath}
\end{equation}
\begin{equation}
    \sum_{t \in [\mathcal{T}]}\sum_{\alpha \in \delta^-((p,t))}x_{\alpha}^{b}
    \leq 1 \quad \forall b \in B \text{, }p \in P \label{constr:pathElementarity}
\end{equation}
\begin{equation}
    f_{a = \pi_\mathcal{D}(\alpha,o)}^{m}= q_{m}x_{\alpha}^{b}\quad \forall b \in
    B \text{, }o \in b\text{, }m \in o\text{, }\alpha \in \mathscr{A}\label{constr:projTTonTS}
\end{equation}
The elementarity constraint is required, as relaxing it allows paths that revisit the same node at different time steps, thereby artificially enhancing consolidation opportunities.

\paragraph{Transportation and Bin-Packing}
Transport on arcs $a \in \mathcal{A}$ involves consolidation.
The shipper procures transport capabilities, referred to generically as \emph{services}, which typically correspond to trucks or shipping containers.
Those services exactly corresponds to the \emph{bins} in a bin-packing problem.
All commodities assigned to a given arc must be packed into these bins, requiring the solution of a bin-packing problem.
In practice, this problem is multi-dimensional and subject to additional constraints.
However, for tractability, we approximate it using a classical single-dimensional bin-packing formulation.
We assume that each bin is fully loaded at the origin and completely unloaded at the destination of its respective arc.

For each arc $a \in \mathcal{A}$, let $L_{a} \in \mathbb{R}$ denote the bin capacity and $K_a \in \mathbb{N}$ the maximum number of bins available.
In our context, this constant is large enough in order to fulfill all the services needs of an instance, as an automotive inbound supply chain is not limited by a single carrier fleet capacity.
We introduce decision variables $\tau_{a}^{k} \in \{0, 1\}$ to indicate whether bin $k \in [K_a]$ is used, and $y_{ak}^{m} \in \mathbb{N}$ to represent the number of units of commodity $m$ assigned to service $k$.
This bin-packing model yields the constraints~\eqref{constr:binDispatch} and~\eqref{constr:binCapa}.
\begin{equation}
    f_{a}^{m}= \sum_{k \in K}y_{ak}^{m}\quad \forall m \in M \text{, }a \in
    \mathcal{A}^{\mathrm{con}}\label{constr:binDispatch}
\end{equation}
\begin{equation}
    \sum_{m \in M}y_{ak}^{m}\ell_{m}\leq L_{a}\tau_{a}^{k}\quad \forall a \in
    \mathcal{A}^{\mathrm{con}}\text{, }k \in [K_a]\label{constr:binCapa}
\end{equation}

\paragraph{Outsourcing Exception}
Transportation on a subset of arcs, specifically $a \in \mathcal{A}^{\mathrm{out}} \subset \mathcal{A}$, is \emph{outsourced}.
In this case, the shipper procures transport services from a third-party logistics provider (3PL) on a per-commodity basis, rather than contracting directly with a carrier.
As a result, there is no need to explicitly model bin-level consolidation on these arcs.
To distinguish between outsourced and consolidated flows, we define the set of \emph{consolidated arcs} $\mathcal{A}^{\mathrm{con}} \subseteq \mathcal{A}$, which comprises all arcs requiring explicit bin-packing and consolidation.
We therefore have $\mathcal{A}^{\mathrm{out}} \;\uplus\; \mathcal{A}^{\mathrm{con}} = \mathcal{A}$.
In Renault's case, we even have $\mathcal{A}^{\mathrm{out}} \subseteq \mathcal{A}^{\mathrm{col}}$.

\subsection{Cost Structure and Full Problem Formulation}\label{subsec:cost_and_full_pb}

This subsection presents the complete cost structure and optimization formulation for the shipper transportation planning problem.

\paragraph{Transportation}

The model assigns a \emph{transport cost} for each bin used on consolidated arcs $a \in A^{\mathrm{con}}$, denoted by $c_{a}^{\mathrm{con}}$.
On outsourced arcs $a \notin A^{\mathrm{con}}$, the shipper incurs a volume-based cost $c_{a}^{\mathrm{out}}$.
Each arc also incurs a volume-dependent cost $c_{a}^{\mathrm{CO2}}$ based on its carbon emission factor.

\paragraph{Platforms}

Processing and handling goods at logistics platforms generates costs.
The platform cost $c_{p}^{\mathrm{plat}}$ increases proportionally with the volume processed.
Platform contracts define a capacity limit $u_{p}^{t}$ for each period and $u_{p}^{[T]}$ for the whole planning horizon, specifying the maximum volume allowed at the base rate.
Exceeding this capacity leads to an overload cost $c_{p}^{\mathrm{over}}$.
This representation best reflects Renault’s contractual and operational practices, in accordance with logistics and platform‑management teams.
Let $z_{p}^{t} \in \mathbb{R}$ represent the excess volume at platform $p$ for each period and $z_p^{[T]}$ the excess over the horizon.
The model enforces the following overload constraints:
\begin{equation}
    \sum_{a \in \delta^-(p)}\sum_{m \in M}f_{a}^{m}\ell_{m} \leq u_{p}^{t} + z_{p}^{t}
    \quad \forall p \in \mathcal{P}
    \label{constr:platOverloadStep}
\end{equation}
\begin{equation}
    \sum_{a \in \mathcal{A}(\delta^-(p))}\sum_{m \in M}f_{a}^{m}\ell_{m} \leq u_{p}^{[T]} + z_{p}^{[T]}
    \quad \forall p \in P
    \label{constr:platOverloadHorizon}
\end{equation}
Constraint~\ref{constr:platOverloadHorizon} represents a contractual volume threshold negotiated with platform operators over the planning horizon. 
It primarily serves to model overload tariffs and improve procurement-cost estimation.

\paragraph{Commodities}

Commodities in transit tie up capital, which we model using a \emph{capital cost} $c_{m}^{\mathrm{cap}}$, proportional to the distance traveled.
The cost incurred by commodity $m$ on arc $a$ is given by $c_{am}^{\mathrm{cap}} = d_{a} \times c_{m}^{\mathrm{cap}}$.

\paragraph{Network Cost Function}

To simplify notation, we aggregate all commodity-specific costs for arc $a = (u,v)$ into a composite term $c_{am}^{\mathrm{com}} = \frac{\ell_{m}}{L_{a}}c_{a}^{\mathrm{CO2}} + c_{am}^{\mathrm{cap}} + \mathbb{I}(v \in \mathcal{P}) \frac{\ell_{m}}{L_{a}}c_{v}^{\mathrm{plat}} + \mathbb{I}(a \notin \mathcal{A}^{\mathrm{con}}) \frac{\ell_{m}}{L_{a}}c_{a}^{\mathrm{out}}$.
Here, $\mathbb{I}(v \in \mathcal{P})$ equals 1 if node $v$ is a platform, and 0 otherwise.
Likewise, $\mathbb{I}(a \notin \mathcal{A}^{\mathrm{con}})$ equals 1 if arc $a$ is outsourced.
We also aggregate platform overloading volumes into a single variable $z_p^\mathrm{total} = \sum_{\nu \in \mathcal{V}(p)} z_{\nu}^{t} + z_{p}^{[T]}$.
Then, the total network cost is:
\begin{equation}
    \sum_{a \in \mathcal{A}} \left( \sum_{k \in K} c_{a}^{\mathrm{con}} \tau_{a}^{k}
    + \sum_{m \in M} c_{am}^{\mathrm{com}} f_{a}^{m} \right)
    + \sum_{p \in P} c_{p}^{\mathrm{over}} z_{p}^{\mathrm{total}}
    \label{obj:fullNetworkCost}
\end{equation}
This objective captures the cost of using bins on consolidated arcs, volume-dependent costs on all arcs, and overload penalties at platforms.

\paragraph{Shipper Transportation Planning Problem}

The complete \emph{Shipper Transportation Planning} Problem can now be formulated as the following mixed-integer program:
\begin{equation}
    \label{pb:forwardPath}
    \begin{aligned}
        \min_{x,y,z,f,\tau} \quad & \text{Network Cost}~\eqref{obj:fullNetworkCost} \\
        \text{s.t.} \quad 
        & \text{Flow Constraints}~\eqref{constr:bundlePath}, ~\eqref{constr:pathElementarity},~\eqref{constr:projTTonTS} \\
        & \text{Packing Constraints}~\eqref{constr:binDispatch},~\eqref{constr:binCapa} \\
        & \text{Overload Constraints}~\eqref{constr:platOverloadStep},~\eqref{constr:platOverloadHorizon} \\
        & \text{Flow variables } x \in \{0,1\},\ f \in \mathbb{N} \\
        & \text{Packing variables } \tau \in \{0,1\},\ y \in \mathbb{N} \\
        & \text{Platform variables } z \in \mathbb{R}_{+}
    \end{aligned}
\end{equation}

Even with size-reduction techniques, solving industrial-scale instances remains intractable.
Typical instances feature approximately $|\mathcal{A}^{\mathrm{con}}||K||M| \approx 10^{12}$ variables and $|B||\mathscr{V}| \approx 10^{10}$ constraints, far exceeding the capability of current (commercial) solvers.

\paragraph{Discussion and Extensions}

Although tailored to Renault’s operations, the model accommodates a wide range of extensions, including multi-modal transport, vehicle tours, inventory, part sourcing and return logistics.
Appendix~\ref{appendix:model_extensions} provides details for implementing these extensions.
While the remainder of this paper focuses on Renault’s specific operational setting, 
the proposed model and the solution approach developed accommodates with ease the additional features, enabling their usage to a broad range of planning environments.
    
\section{Decomposition-guided ILS} \label{sec:ILS}

The problem instances under consideration involve far too many variables and constraints,  
rendering exact methods such as mixed-integer programming or classical decomposition techniques computationally infeasible.
To address this challenge, we develop a tailored \emph{Iterated Local Search} (ILS) algorithm designed for large-scale combinatorial optimization under structural constraints.

\begin{figure}[b]
    \centering
    \includegraphics[width=\textwidth]{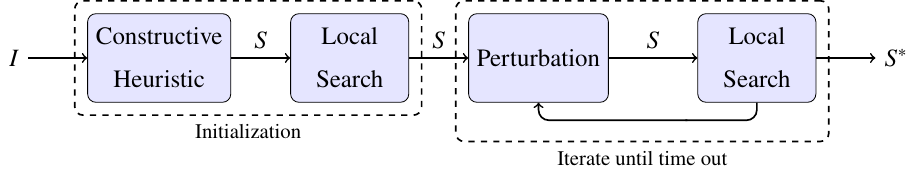}
    \caption{Description of the Iterated Local Search} 
    \label{fig:ILS}
    \vspace{-0.25cm}
    \note{A constructive heuristic and a first local search finds a promising candidate solution.
    The search space is then explored and the solution refined by alternating between a local search and a perturbation phase.}
\end{figure}

Described in Figure~\ref{fig:ILS}, it works
as follows.
Given an instance $I$, a \emph{constructive} heuristic first builds an initial
solution $S$, which is then refined by alternating between a \emph{local search}
that improves the current solution locally and a \emph{perturbation} phase that gets out of local minima. 
As discussed in the case study and literature review, these problems differ fundamentally from the carrier-centric formulations that dominate the service network design literature: these problems procure rather than operate transportation resources, focus on long-term transportation policies, have larger and sparser networks and enables detailed consolidation modeling and flexible delivery schedules.
These structural properties strongly correlates with four key observations: (i) the problem becomes tractable when restricted to a single commodity bundle; (ii) Due to consolidation effects, bundles rarely follow their individual shortest paths in the network; (iii) A significant portion of transportation costs is concentrated on shared network segments; (iv) Bundles relying on these segments are typically intercontinental or low-volume long-distance, making them especially dependent on consolidation.
The first observation is derived from a mathematical analysis of the problem while the remaining ones originate from the analysis of Renault's operational transportation data.

Based on these insights, we decompose the overall algorithm and its components into three main types of computational routines: packing, single-bundle, and multi-bundle operations.
Each component, defined as a variant of a routine, plays a specific role within the ILS framework—constructing initial solutions, refining existing ones, or exploring larger neighborhoods via controlled perturbation.
Although our algorithm builds on a standard Iterated Local Search framework, its performance derives not from the metaheuristic structure itself but from operators and perturbation mechanisms specifically designed for the geometric and structural features of this problem. 
The constructive and reinsertion operators explicitly exploit the dual time-expanded structure to transform complex bundle-routing decisions into shortest-path problems. 
The consolidate and refine neighborhoods are specifically designed to uncover alternative consolidation and routing patterns.
To address the combinatorial interactions between routing decisions, path-based perturbation mechanisms deliberately redirect bundles through common portions of the network.
Our contribution therefore lies in two key design principles: it systematically exploits the tractability of bundle-level routing subproblems to explore large neighborhoods at low computational cost, and it explicitly targets the shared transportation corridors that concentrate consolidation opportunities, enabling the metaheuristic to navigate the solution space effectively despite the problem’s scale and combinatorial complexity.
Table~\ref{tab:ils_blocks} summarizes the computational building blocks and their respective variants, including the geometric principles they exploit and their functional roles within the algorithm.

\begin{table}[!t]
\centering
\caption{Computation types and their variants.} \label{tab:ils_blocks}
\renewcommand{\arraystretch}{1.3} 
\resizebox{\textwidth}{!}{%
\begin{tabular}{>{\raggedright\arraybackslash}p{3cm} | >{\raggedright\arraybackslash}p{3cm} | >{\raggedright\arraybackslash}p{4.5cm} | >{\centering\arraybackslash}p{2cm} | >{\raggedright\arraybackslash}p{5.5cm}}
\toprule
\textbf{Type} & \textbf{Variants} & \textbf{Geometric Insight} & \textbf{Usage} & \textbf{Description} \\
\midrule
\multirow{5}{*}{Packing} 
  & Pack & Batched Bin-Packing & All 
  & Inserts commodities on arcs; computes associated cost \\
  \cmidrule{2-5}
  & Re-Pack & Classical Bin-Packing & LS 
  & Removes and re-packs all commodities to allow cost re-optimization \\
\midrule
\multirow{4}{*}{Single Bundle} 
  & Insertion & Shortest Path on Sparse Graph (weighted with Pack) & CH 
  & Builds paths one commodity at a time using sparse cost structure \\
  \cmidrule{2-5}
  & Re-Insertion & Shortest Path on Sparse Graph (weighted with Pack) & LS 
  & Removes and reinserts bundles to escape local minima \\  
\midrule
\multirow{6}{*}{Multi-Bundle} 
  & Consolidate & Shortest Path on Sparse Graph (weighted with Pack) & LS 
  & Merges bundles on shared paths before re-insertion; partial path recomputation \\
  \cmidrule{2-5}
  & Flow-based & Partial MILP Relaxation & P 
  & Solves aggregated routing via relaxed MILP to explore new configurations \\
  \cmidrule{2-5}
  & Path-based & Partial MILP Relaxation & P 
  & Similar to flow-based, but works on enumerated route alternatives \\
\bottomrule
\end{tabular}
}
\vspace{-0.25cm}
\note{Usage: CH = Constructive Heuristic, LS = Local Search, P = Perturbation.}
\end{table}

\subsection{Constructive algorithm}\label{subsec:constructive}

\begin{figure}[!t]
    \centering
    \includegraphics[width=0.95\textwidth]{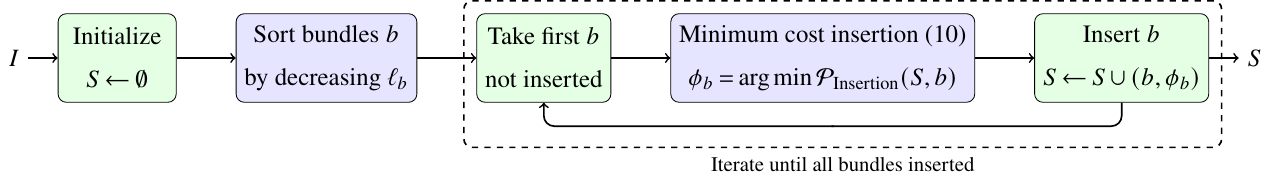}
    \caption{Description of the Constructive Heuristic} 
    \label{fig:constructive}
    \vspace{-0.25cm}
    \note{Bundles are sorted by decreasing size and then inserted one by one in the current solution.}
\end{figure}

The constructive heuristic incrementally builds a feasible solution \( S \) for a given instance \( I \) by inserting commodities bundle by bundle.
At each step, it chooses configurations that minimize arc activation and platform usage costs by solving a bundle insertion sub-problem.
Figure~\ref{fig:constructive} illustrates the heuristic's logic.
The heuristic follows the classical \emph{First-Fit Decreasing} (FFD) rule from bin-packing~\citep{korf2002ffd}.
FFD processes items in decreasing order of size and inserts each into the first bin with sufficient space.
In our setting, the heuristic inserts bundles \( b \in B \) into the partial solution \( S \) in decreasing order of their maximum packaging size \( \ell_B = \max_{m \in B} \ell_m \), which serves as a proxy for packing complexity.
At each iteration, the heuristic inserts a bundle in the current solution using a minimum cost insertion subproblem and maintains a partial solution that satisfies constraints~\eqref{constr:bundlePath}--\eqref{constr:binCapa} for the current set of bundles \( B^S \subset B \).
Each bundle \( b \in B^S \) is associated with a path \( \phi_b \), and the set of inserted commodities grows as the algorithm progresses.

\paragraph{Minimum Cost Insertion}

To insert a new bundle \( b \) into the partial solution \( S \) at minimum cost, the algorithm solves a restricted version of problem~\eqref{pb:forwardPath}.
This subproblem identifies a feasible delivery path for \( b \) and evaluates the resulting consolidated arc costs as follows: define \( M^S = \bigcup_{b \in B^S} \{m \in b\} \) as the set of commodities already in \( S \), and let \( (\bar{y}, \bar{f}) \) denote their fixed packing decisions. 
Using those, it enforces solution-aware packing and platform constraints when inserting~\( b \): bin capacity $L_a$ becomes $\displaystyle \bar{L}_a^k = L_a - \sum_{m \in M^S} \bar{y}_{ak}^m \ell_m$ and platform capacity $u_p^{t}$ and $u_p^{[T]}$ becomes $\displaystyle \bar{u_p}^{t} = u_p^{t} - \sum_{a \in \delta^-(p)}\sum_{m \in M^S} f_a^m \ell_m$ for all $p \in \mathcal{P}$ and $\displaystyle \bar{u_p}^{[T]} = u_p^{[T]} - \sum_{a \in \mathcal{A}(\delta^-(p))}\sum_{m \in M^S} f_a^m \ell_m$ for all $p \in P$.
The insertion subproblem $\mathcal{P}_\mathrm{Insertion}(S,b)$ is defined as follows:
\begin{equation}
\label{pb:bundlePricing}
\begin{aligned}
    \min \quad & \text{Network Cost}~\eqref{obj:fullNetworkCost} \\
    \text{s.t.} \quad
    & \text{Flow constraints}~\eqref{constr:bundlePath},~\eqref{constr:pathElementarity},~\eqref{constr:projTTonTS} \text{ for } b \\
    & \text{Packing constraints}~\eqref{constr:binDispatch},~\eqref{constr:binCapa} \text{ with } \bar{L}_a^k \\
    & \text{Platform constraints}~\eqref{constr:platOverloadStep},~\eqref{constr:platOverloadHorizon} \text{ with } (\bar{u}_p^t, \bar{u}_p^{[T]}) \\
    & x, \tau \in \{0,1\}, \quad y, f \in \mathbb{N}, \quad z \in \mathbb{R}_+
\end{aligned}
\end{equation}

To simplify notation, the algorithm assumes all decision variables for bundles not yet inserted are set to zero.

\paragraph{Solving the Minimum Cost Insertion}

The algorithm solves problem~\eqref{pb:bundlePricing} efficiently by reducing it to a shortest path problem on a specialized subgraph of \( \mathscr{G} \).
Since the subproblem considers a single bundle, the algorithm can precompute packing costs and project them from \( \mathcal{G} \) onto \( \mathscr{G} \).
The bundle-specific digraph \( \mathscr{G}^b = (\mathscr{V}^b, \mathscr{A}^b) \) contains only nodes and arcs relevant to bundle \( b \), defined as $\mathscr{V}^b = \{ \nu = (v,t) \in \mathscr{V} \mid v \in \{s_b, u_b\} \cup P \}, \;
\mathscr{A}^b = \{ \alpha = (u,v) \in \mathscr{A} \mid u, v \in \mathscr{V}^b \}$
In practice, the set $\mathscr{V}^b$ is further pruned by removing platforms that cannot be reached from $s_b$ or cannot reach $u_b$.
For each arc \( \alpha \in \mathscr{A}^b \), the algorithm computes three types of costs: 
the \emph{commodity cost} $\displaystyle c_{\alpha}^{\mathrm{com}} = \sum_{\substack{o \in b \\ a = \pi_{\mathcal{D}} (\alpha,o)}}
\ell_o \, c_{am}^{\mathrm{com}} 
$,
the \emph{platform overloading cost} $ \displaystyle 
c_{\alpha}^{\mathrm{over}} = c_p^{\mathrm{over}} \left [ \ell_b - u_p^{[T]} \right ]_+ + \sum_{\substack{o \in b \\ a = \pi_{\mathcal{D}} (\alpha,o)}}
c_p^{\mathrm{over}} \left [ \ell_o - u_p^{t} \right ]_+ ,
$
if arc \( \alpha \) ends at a platform \( p \in P \) and finally the \emph{consolidated arc cost}
\( \displaystyle 
c_{\alpha}^{\mathrm{con}} = \sum_{\substack{o \in b \\ a = \pi_{\mathcal{D}}(\alpha, o)}}
c_{a,o}^{\mathrm{con}}
\).
The consolidated cost is computed with a bin-packing for each projected arc.
\[
c_{a,o}^{\mathrm{con}} = 
\left\{
\begin{aligned}
    \min \quad & \sum_{k \in K} c_a \tau_a^k \\
    \text{s.t.} \quad & q_m = \sum_{k \in K} y_{ak}^m && \forall m \in o \\
    & \sum_{m \in o} y_{ak}^m \ell_m \leq \bar{L}_a^k \tau_a^k && \forall k \in [K_a] \\
    & \tau \in \{0,1\}, \quad y \in \mathbb{N}
\end{aligned}
\right.
\]
Each arc \( \alpha \in \mathscr{A}^b \) is then weighted by the sum if these three components
\(
    c_{\alpha}= c_{\alpha}^{\mathrm{com}}+ c_{\alpha}^{\mathrm{over}}+ c_{\alpha}
    ^{con}
\)

\begin{prop}
    Solving the minimum cost insertion problem $\mathcal{P}_\mathrm{Insertion}(b)$~\eqref{pb:bundlePricing}
    amounts to computing an elementary shortest path from
    $(s_{b}, \mathcal{T}- \tau_{b})$ to $(u_{b}, \mathcal{T})$ in
    $\mathscr{G}^{b}$.
\end{prop}

In practice, this result leads to two important implementation choices.
First, the elementarity constraint~\eqref{constr:pathElementarity} is nearly always inactive.
In over 99\% of cases, the bundle path is already elementary, allowing the use of Dijkstra’s algorithm with non-negative arc costs to efficiently solve the problem. 
Second, in the rare cases where elementarity must be enforced, the algorithm insert the bundle in an empty solution instead of the current solution, which guarantees that the path is going to be elementary.

To maintain tractability across the large number of arc-level bin-packing problems—up to \( O(10^5) \) instances—the algorithm applies the \textsc{FirstFitDecreasing} heuristic.
This approach ensures near-optimal performance in practice while keeping computational effort manageable.
It also satisfies a known approximation guarantee~\citep{Dsa2007TheTB}: 
$ 
\forall I, \ \text{FFD}(I) \leq \frac{11}{9} \cdot \text{OPT}(I) + \frac{6}{9}.
$

\subsection{Local Search}

\begin{figure}[!b]
    \centering
    \includegraphics[width=0.95\textwidth]{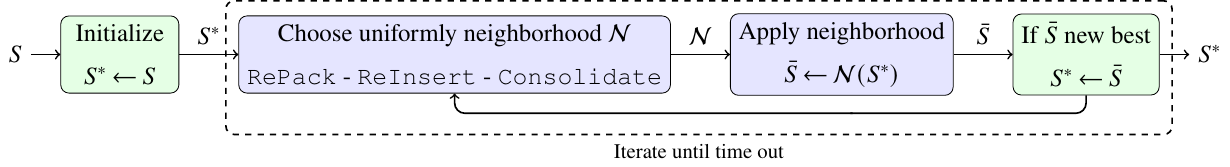}
    \caption{Description of the Local Search procedure} 
    \label{fig:localsearch}
    \vspace{-0.25cm}
    \note{Neighborhoods are uniformly sampled to produce candidate neighbors, becoming the current solution in case of improvement.}
\end{figure}

The local search algorithm operates through three custom-designed and increasingly large neighborhoods: \emph{Re-Pack}, \emph{Re-Insert}, and \emph{Consolidate-and-Refine}.
At each iteration, it randomly selects one of these neighborhoods and generates a neighboring solution.
If this neighbor improves upon the current solution, the algorithm accepts it as the new incumbent; otherwise, it discards the candidate and continues the search.
Figure~\ref{fig:localsearch} illustrates this procedure.

The neighborhood designs draw on empirical observations from our numerical experiments.
While local search methods often benefit from well-crafted neighborhoods, \citet{turkes2021} show that exhaustive tuning of neighborhood mechanics rarely yields significant gains.
We therefore prioritized simplicity and computational efficiency.
Implementation optimizations involves mostly the parallelism of arc cost computations and the efficient allocation of shared memory for bin-packing computations.
A detailed discussion on computational efficiency considerations is available in Appendix~\ref{appendix:computational_efficiency}.
Instead of exhaustively evaluating all neighbors, the algorithm samples a single random neighbor in each iteration.

\begin{figure}[!b]
    \centering
    \includegraphics[width=0.95\textwidth]{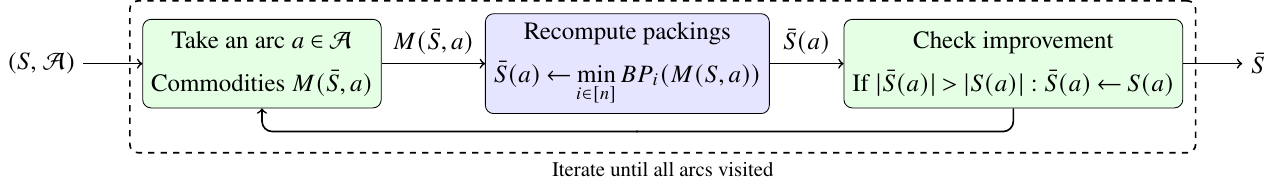}
    \caption{Description of the Re-Pack Neighborhood} 
    \label{fig:re_pack}
    \vspace{-0.25cm}
    \note{Given a set of arcs $\mathcal{A}$, new bins are computed based on provided bin-packing heuristics.}
\end{figure}

\paragraph{Re-Pack Neighborhood}
In this neighborhood, the algorithm recomputes bin-packings across a set of arcs $\mathcal{A}$ given using a portfolio of heuristics \( (BP_i)_{i \in [n]} \).
In practice, two heuristics are used: \emph{First-Fit Decreasing} and \emph{Best-Fit Decreasing}.
By globally recomputing packings, the algorithm avoids suboptimal local packing sequences that arise from purely sequential insertions.
Figure~\ref{fig:re_pack} visualizes this mechanism.
Since this operation cannot increase the cost of the current solution, it is always accepted when feasible.

\begin{figure}[!t]
    \centering
    \includegraphics[width=0.95\textwidth]{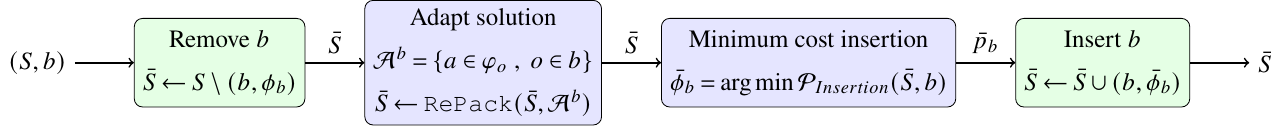}
    \caption{Description of the Re-Insert Neighborhood} 
    \label{fig:re_insert}
    \vspace{-0.25cm}
    \note{Building on the constructive heuristic, a bundle is removed and re-inserted in the adapted solution.}
\end{figure}

\paragraph{Re-Insert Neighborhood}
This neighborhood removes a randomly selected bundle from the current solution, adapts the affected bins and reinserts it using the insertion procedure \( \mathcal{P}_{Insertion}(S, b) \) defined in~\eqref{pb:bundlePricing}.
Figure~\ref{fig:re_insert} shows the process. 

\begin{figure}[!b]
    \centering
    \includegraphics[width=0.95\textwidth]{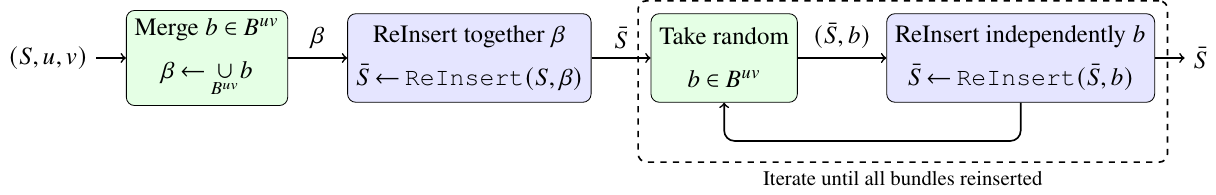}
    \caption{Description of the Consolidate-and-Refine Neighborhood} 
    \label{fig:consolidate}
    \vspace{-0.25cm}
    \note{This neighborhood seeks consolidation for bundles flowing from $u$ to $v$ on this part of their paths.}
\end{figure}

\paragraph{Consolidate-and-Refine}
The third neighborhood, illustrated in Figure~\ref{fig:consolidate}, targets consolidation commodities along the common network connecting platforms and industrial sites.
It selects two nodes \( u \) and \( v \) and collects the set of bundles \( B^{uv} \) that travel between them.
It then aggregates these bundles into a temporary super-bundle \( \beta = \bigcup_{b \in B^{uv}} b \) and inserts it as a single entity using the \texttt{ReInsert} operator.
After this coarse consolidation, the algorithm re-evaluates each bundle in \( B^{uv} \) individually via further \texttt{ReInsert} operations to refine the solution.
This ensures that only bundles benefiting from consolidation remain on the new path.

\subsection{Perturbation Scheme} \label{subsubsec:perturb}

\begin{figure}[!b]
    \centering
    \includegraphics[width=0.9\textwidth]{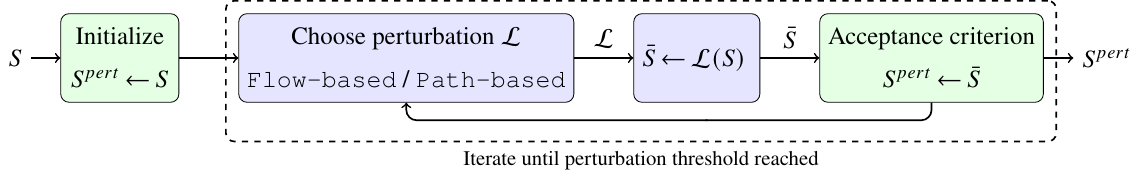}
    \caption{Description of the Perturbation Scheme} 
    \label{fig:perturbation}
    \vspace{-0.25cm}
    \note{Flow-based and Path-based perturbations are alternated until the perturbation threshold is reached.}
\end{figure}

The final component of our Iterated Local Search metaheuristic is a perturbation scheme.
This module acts as a large-neighborhood search mechanism based on integer programming.
It enables the algorithm to escape local minima by re-optimizing parts of the solution and exploring promising regions of the search space.
We refer to these large-scale modifications as \emph{perturbations} because they rely on approximate subproblems.
Although this approximation increases tractability, it also means that solutions generated through perturbation may not always improve upon the incumbent.
A cost increase tolerance is therefore set 1.5\% to discard perturbations that would degrade the solution quality too much.
Figure~\eqref{fig:perturbation} illustrates the perturbation process.

\paragraph{Perturbation Design}

Each perturbation re-optimizes a subset of bundles \( \bar{B} \subseteq B \), following a principle similar to the \texttt{ReInsert} operator, but on a larger scale.
The algorithm removes all bundles in \( \bar{B} \) from the current solution and fixes the configuration of all other bundles in \( B \setminus \bar{B} \).
It then solves a simplified MILP to reinsert \( \bar{B} \).
This simplification replaces the detailed packing constraints~\eqref{constr:binDispatch} and~\eqref{constr:binCapa} with a giant container approximation:
\begin{equation}
\bar{L}_{a} + \sum_{m \in \bar{B}} f_{a}^{m} \ell_{m} \leq L_{a} \tau_{a}, \quad \tau_{a} \in \mathbb{N}, \quad \forall a \in \mathcal{A}^{\mathrm{con}}
\label{constr:binCapaGiant}
\end{equation}
To mitigate the effect of this relaxation, we scale the transport cost on arcs $a \in \mathcal{A}^{\mathrm{con}}$ using a slope-scaling heuristic inspired by \citet{jarrah2009}.
The adjusted arc cost becomes 
\( \displaystyle
\bar{c}_{a} = c_{a} \cdot \frac{BP(a)}{\lceil \sum_{m \in M} f_{a}^{m} \ell_{m} \rceil}
\).
The resulting objective function is given by:
\begin{equation}
\sum_{a \in \mathcal{A}} \left( \bar{c}_{a} \tau_{a} + \sum_{m \in M} c_{am}^{\mathrm{com}} f_{a}^{m} \right) + \sum_{p \in \mathcal{P}} c_{p}^{\mathrm{over}} z_{p}
\label{obj:approxNetCost}
\end{equation}

We consider two formulations for implementing perturbations: \emph{flow-based reoptimization} and \emph{path-based reoptimization}.
The former allows complete routing flexibility but is tractable only for small \( \bar{B} \); the latter restricts routing to a small number of preselected paths per bundle, allowing it to scale to larger \( \bar{B} \).

\paragraph{Flow-Based Reoptimization}

This approach reinserts a subset of bundles \( \bar{B} \) by solving the following MILP:
\begin{equation}
\label{pb:perturbMilp}
\begin{aligned}
\min \quad & \text{Approximate Network Cost}~\eqref{obj:approxNetCost} \\
\text{s.t.} \quad
& \text{Flow constraints}~\eqref{constr:bundlePath},~\eqref{constr:pathElementarity},~\eqref{constr:projTTonTS} \text{ for } \bar{B}, \\
& \text{Packing constraints}~\eqref{constr:binCapaGiant} \text{ for } (S, \bar{B}), \\
& \text{Platform constraints}~\eqref{constr:platOverloadStep},~\eqref{constr:platOverloadHorizon} \text{ for } (S, \bar{B}), \\
& x \in \{0,1\}, \quad f, \tau \in \mathbb{N}, \quad z \in \mathbb{R}_+
\end{aligned}
\end{equation}

We define three families of perturbations within this framework, each based on a shared bundle characteristic (Table~\ref{tab:milpArcSubProblems}).
To maximize the effectiveness of this scheme, the algorithm stacks multiple subproblems until the aggregated instance reaches a variable budget that can be solved within a few minutes (on the order of several million variables).
Perturbations continue iteratively until the number of modified paths exceeds a specified threshold.

\paragraph{Path-Based Reoptimization}

This alternative formulation improves scalability for larger \( \bar{B} \).
Instead of computing paths during optimization, it assumes a small set of feasible paths \( \varphi \in \Phi^b \) is given for each bundle \( b \in \bar{B} \).
The flow constraints are replaced by path-based formulations:
\begin{equation}
\sum_{\varphi \in \Phi^b} x_{\varphi}^b = 1, \quad x_{\varphi}^b \in \{0,1\}, \quad \forall b \in \bar{B}
\label{constr:oneBundlePath}
\end{equation}
\begin{equation}
f_{\pi_{\mathcal{D}}(\alpha, o)}^m = q_m \sum_{\varphi \ni \alpha} x_{\varphi}^b, \quad \forall (b, o, m) \in \bar{B}, \; \alpha \in \mathscr{A}
\label{constr:pathProj}
\end{equation}

We define three families of path-based perturbations (Table~\ref{tab:milpPathSubProblems}), two of them adapted from strategies proposed by \citet{lindsey2016improved}.
The \emph{Attract} perturbation introduces bundles to arcs they currently avoid, while \emph{Reduce} attempts to reroute bundles away from congested arcs.
Finally, the \emph{Directs} perturbation is used as the last perturbation before activating the local search as it tends to switch bundles on directs paths to the shared network.

\begin{table}[!t]
\centering
\caption{Three families of flow-based reoptimization.} \label{tab:milpArcSubProblems}
\renewcommand{\arraystretch}{1.3} 
\begin{tabular}{>{\raggedright\arraybackslash}p{3cm} | >{\raggedright\arraybackslash}p{4cm} | >{\raggedright\arraybackslash}p{6cm}}
\toprule
\textbf{Name} & \textbf{Family of Subproblem} & \textbf{Subproblem \( \bar{B} \)} \\
\midrule
Single Plant & Industrial site \( u \in U \) & \( \{ b \in B \mid u_b = u \} \) \\
\midrule
Single Supplier & Supplier \( s \in S \) & \( \{ b \in B \mid s_b = s \} \) \\  
\midrule
Random & Number \( n \geq 1 \) & Uniform sample of size \( n \) from \( B \) \\
\bottomrule
\end{tabular}
\end{table}

\begin{table}[!t]
\centering
\caption{Three families of path-based reoptimization.} \label{tab:milpPathSubProblems}
\renewcommand{\arraystretch}{1.3} 
\begin{tabular}{>{\raggedright\arraybackslash}p{3cm} | >{\raggedright\arraybackslash}p{4cm} | >{\raggedright\arraybackslash}p{5cm}}
\toprule
\textbf{Name} & \textbf{Family of Subproblem} & \textbf{Subproblem \( \bar{B} \)} \\
\midrule
Attract & Arc \( a \in \mathscr{A}^{\mathrm{pla}} \cup \mathscr{A}^{\mathrm{del}} \) & \( \{ b \in B \mid a \notin p_b \} \) \\
\midrule
Reduce & Arc \( a \in \mathscr{A}^{\mathrm{pla}} \cup \mathscr{A}^{\mathrm{del}} \) & \( \{ b \in B \mid a \in p_b \} \) \\  
\midrule
Directs & $\emptyset$ & \( \{ b \in B \mid b \text{ on direct path} \} \) \\
\bottomrule
\end{tabular}
\end{table}

As in the flow-based case, we stack multiple subproblems to scale the perturbation.
Since not all bundles are eligible for rerouting through a given arc, we prioritize arcs where a large share of bundles are \emph{candidates}—i.e., bundles for which the arc is either currently unused (\emph{Attract}) or already used (\emph{Reduce}).
While these MILPs are typically faster to solve, they require extensive path enumeration, which can offset their computational advantages at the scale of a full perturbation.

\section{Tractable Lower Bounds}\label{sec:lowerbound}
\begin{figure}[!b]
    \centering
    \includegraphics[scale=1.0]{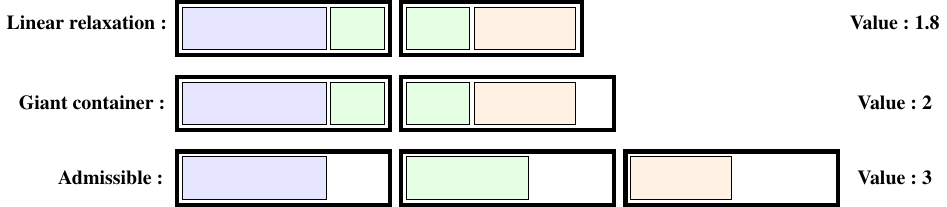}
    \caption{Example of linear relaxation (top) and giant container relaxation (middle) for bin-packing (bottom).} 
    \label{fig:lower-bound-sketch}
\end{figure}

To evaluate the quality of the proposed heuristic, we consider three lower bounding procedures that offer different trade-offs between tightness and computational complexity.
These procedures, summarized in Figure~\ref{fig:lower-bound-sketch}, correspond to successive relaxations of the bin-packing structure: a linear relaxation, a mixed giant container relaxation, and a full giant container relaxation.
The linear relaxation can be decomposed into independent subproblems per bundle, making it computationally tractable but relatively weak.
The mixed giant container relaxation strengthens this bound on direct arcs, offering a more favorable balance between quality and efficiency.
Extending the giant container relaxation to all arcs yields the full giant container bound, which yields a stronger bound at the expense of increased computational costs.
Finally, from the mixed giant container relaxation we derive a rounding heuristic that produces feasible solutions and serves as a benchmark for assessing the performance of our ILS algorithm.

\paragraph{Linear Relaxation Bound.}

We first consider the linear relaxation of problem~\eqref{pb:forwardPath}. 
In this relaxation, the packing variables $\tau$ and $y$ as well as the packing constraints~\eqref{constr:binCapa} and~\eqref{constr:binDispatch} can be removed.  
These elements serve only to determine the number of transport services used by commodities on each arc. 
In the relaxed setting, this number equals the total commodity volume on the arc divided by the service capacity. 
We therefore obtain the following linear relaxation:
\begin{equation}
    \label{pb:linearLB+}
    \begin{aligned}
        \min \, & \sum_{a \in \mathcal{A}} \sum_{m \in M}(c_{am}^{\mathrm{com}}+ \frac{\ell_{m}}{L_{a}}c
        _{a}^{\mathrm{con}}) f_{a}^{m} + \sum_{p \in \mathcal{P}}c_{p}
        ^{\mathrm{over}}z_{p}\\
        \mathrm{s.t.}\, 
        & \ \text{Flow Constraints}~\eqref{constr:bundlePath}, ~\eqref{constr:pathElementarity}\text{ and }~\eqref{constr:projTTonTS} \\
        & \ \text{Platform Overloading Constraints}~\eqref{constr:platOverloadStep},~\eqref{constr:platOverloadHorizon}
    \end{aligned}
\end{equation}
Although this formulation still involves $O(10^{10})$ variables and cannot be solved directly with current linear programming technology, it can be decomposed by bundle if the platform overloading variables $z$ and constraints~\eqref{constr:platOverloadStep} and~\eqref{constr:platOverloadHorizon} are omitted. 
This leads to the following formulation:
\begin{equation}
    \label{pb:linearLB}
    \sum_{b \in B}
    \left \{
    \begin{aligned}
        \min \, & \sum_{a \in \mathcal{A}} \sum_{m \in b}(c_{am}^{\mathrm{com}}+ \frac{\ell_{m}}{L_{a}}c
        _{a}^{\mathrm{con}}) f_{a}^{m} \\
        \mathrm{s.t.}\, 
        & \ \text{Flow Constraints}~\eqref{constr:bundlePath}, ~\eqref{constr:pathElementarity}\text{ and }~\eqref{constr:projTTonTS} \text{ for } b 
    \end{aligned}
    \right . 
\end{equation}

\begin{prop}\label{prop:lowerbound}
    The linear relaxation of~\eqref{pb:forwardPath} has the same value as problem~\eqref{pb:linearLB+}.
    Without platform overloading constraints~\eqref{constr:platOverloadStep} and~\eqref{constr:platOverloadHorizon}, it is equal to the sum of bundle specific subproblems~\eqref{pb:linearLB}.
\end{prop}
Each of these bundle-specific subproblems corresponds to the linear relaxation of the minimum-cost insertion problem~\eqref{pb:bundlePricing} without platform overloading constraints. 
The same solution approach therefore applies: the problem reduces to a shortest path computation on the bundle-specific digraph $\mathscr{G}^b$. 
For the proof of Proposition~\ref{prop:lowerbound} we refer to Appendix~\ref{appendix:lower_bound_proof}.

\paragraph{Mixed-Giant Container Bound.}

The quality of the previous bound can be improved by exploiting a structural property of the problem.
On direct arcs, only commodities from a single bundle are present, meaning that there is no coupling across bundles.
For these arcs, we can therefore strengthen the relaxation by applying a giant container approximation.
The arc cost then becomes
\[
    \forall a \in \mathcal{A}^{dir} \ : \ \left \lceil \sum_{m \in b}\frac{f_{a}^{m}\ell_{m}}{L_{a}} \right \rceil c
        _{a}^{\mathrm{con}} + \sum_{m \in b} c_{am}^{\mathrm{com}} f_{a}^{m}
\]
The non-linearity introduced in the objective function can be handled in practice by precomputing the integer number of transport services required by each order and using this value to evaluate the direct arc costs in the travel-time graph.

\paragraph{Rounding Heuristic.}\label{subsec:rounding}

Both the linear relaxation and the mixed giant container bound produce shortest paths for each bundle. 
By fixing these paths and subsequently solving a bin-packing problem for every arc of the time–space graph, we obtain feasible packings and thus a valid solution to problem~\eqref{pb:forwardPath}. 
We apply this approach using the mixed giant container bound, as it is stronger than the linear relaxation. 
The resulting feasible solution serves as a benchmark for our ILS algorithm. 

\paragraph{Full Giant Container Bound.}
Applying the giant container relaxation to all arcs of the time–space graph yields problem~\eqref{pb:perturbMilp} without cost scaling, where the subset of bundles is $\bar B = B$. 
This formulation inspired the perturbation MILP introduced earlier. 
Unlike the previous bounds, the full giant container relaxation is no longer decomposable by bundle, making it tractable only for instances of moderate size. 
It is closely related to the load plan design problem studied in the literature \citep{erera2013improved,lindsey2016improved}.

\section{Computational Study}
\label{sec:ComputeStudy}
To assess the practical relevance of our approach, we conduct a comprehensive computational study on real-world instances from Renault’s inbound logistics network. 
The analysis proceeds in three steps. 
We first describe the experimental design, including the instances, benchmarks, and implementation details (Section~\ref{subsec:exp_design}). 
Next, we compare solution quality and runtime against alternative strategies to evaluate relative performance (Section~\ref{subsec:perf_analysis}). 
Finally, we extend the analysis to highlight managerial insights on consolidation, regularity, and outsourcing, providing actionable guidance for industrial practice (Section~\ref{subsec:ext_analysis}).
In addition to these core components, the appendix provides two complementary analyses. 
First, we report a detailed hyper‑parameter tuning study, time‑profile evaluations, and an ablation study that quantify the contribution of each ILS component. 
Second, we include an extended version of the managerial and practical analyses, offering deeper insight into integrated versus decomposed planning and the solvability of the model with off‑the‑shelf MILP technology.

\subsection{Experimental design} \label{subsec:exp_design}
To evaluate the performance of our ILS metaheuristic, we address two main questions:  
(1) How do the key difficulties identified in the case study translate into the numerical resolution of real industrial instances?  
(2) How does our algorithm perform compared to established resolution strategies?  
Before answering these questions, we introduce the test instances and present the benchmark heuristics used for comparison.
Description of the implementation and hyperparameter choices are provided in Appendix~\ref{appendix:computational_efficiency} and~\ref{appendix:benchmark_adaptation}.

\begin{figure}[!b]
    \begin{subfigure}{0.33\textwidth}
        \centering
        \includegraphics[width=0.95\textwidth]{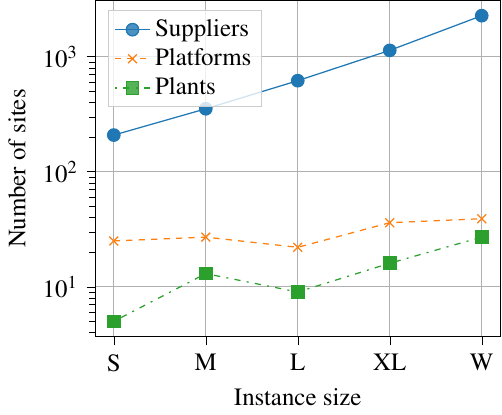}
        \caption{Sites}
    \end{subfigure}
    \begin{subfigure}{0.33\textwidth}
        \centering
        \includegraphics[width=0.95\textwidth]{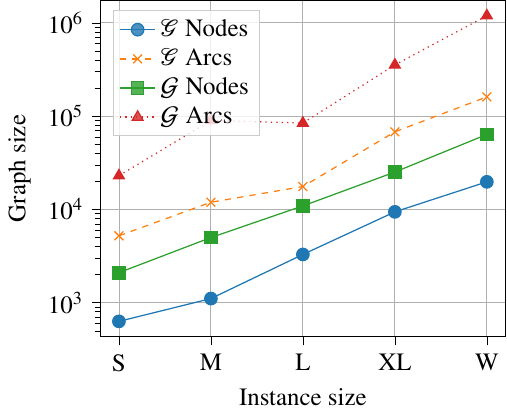}
        \caption{Graph sizes}
    \end{subfigure}
    \begin{subfigure}{0.33\textwidth}
        \centering
        \includegraphics[width=0.95\textwidth]{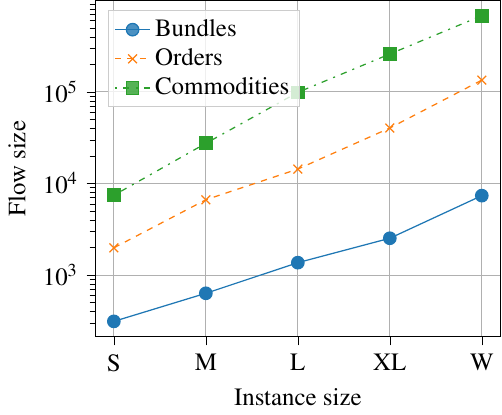}
        \caption{Flows}
    \end{subfigure}
    \caption{Description of scalability test instances.}
    \label{fig:scale-instances}
\end{figure}

\paragraph{Instances}
We rely on two sets of instances derived from real industrial data from Renault's inbound logistics operations.
The first set of instances, illustrated in Figure~\ref{fig:scale-instances}, is designed to address the scalability question.
It consists of five instances of increasing size: \emph{Small (S)}, \emph{Medium (M)}, \emph{Large (L)}, \emph{Very Large (XL)}, and \emph{World (W)}.
The second set of instances, summarized in Table~\ref{tab:perf-instances}, is used for performance evaluation.
It consists of five different \emph{World} instances. 
\begin{table}[!t]
\centering
\caption{Description of performance test instances.} \label{tab:perf-instances}
\renewcommand{\arraystretch}{1.3} 
\resizebox{\textwidth}{!}{\begin{tabular}{>{\raggedright\arraybackslash}p{2cm} | >{\raggedright\arraybackslash}p{1.5cm} >{\raggedright\arraybackslash}p{1.5cm} | >
{\raggedright\arraybackslash}p{1.5cm} >{\raggedright\arraybackslash}p{1.75cm} | >{\raggedright\arraybackslash}p{1.5cm} >{\raggedright\arraybackslash}p{1.5cm} >{\raggedright\arraybackslash}p{2.5cm}}
\toprule
 & \textbf{$\mathscr{G}$ nodes} & \textbf{$\mathscr{G}$ arcs} & \textbf{$\mathcal{G}$ nodes} & \textbf{$\mathcal{G}$ arcs} & \textbf{Bundles} & \textbf{Orders} & \textbf{Commodities} \\
\midrule
Min & 19 689 & 160 466 & 63 518 & 1 192 022 & 7360 & 118 611 & 625 166 \\
\midrule
Max & 19 769 & 161 024 & 63 778 & 1 198 366 & 7473 & 137 506 & 724 454 \\
\bottomrule
\end{tabular}}
\end{table}

\paragraph{Benchmarks} \label{subsec:benchmarks}
We compare the performance and scalability of our ILS algorithm against the following benchmark heuristics:  

\begin{itemize}
    \item \textbf{Shortest (S)} Each bundle is assigned the shortest path in terms of distance.
    \item \textbf{Renault-based (R)} The current transport operations are projected into our modeling framework to produce an admissible solution.   
    \item \textbf{Average (A)} All orders of a bundle are aggregated into a single order delivered at time step~1.
The delivery path is then computed iteratively using the full giant container relaxation.  
    \item \textbf{Constructive (C)} The constructive heuristic described in Section~\ref{subsec:constructive}.  
    \item \textbf{Lower Bound Rounding (LBR)} The rounding heuristic described in Section~\ref{subsec:rounding}.  
    \item \textbf{Slope Scaling (SS)} Adaptation of the heuristic developed by \citet{jarrah2009}, based on iterative linearization of the problem.
    \item \textbf{IP-based Local Search (IPLS)} Adaptation of the heuristic developed by \citet{lindsey2016improved}, based neighbors obtained by exact MILP subproblems reoptimization.
\end{itemize}

The Renault-based solution provides the first natural benchmark. 
As it reflects the industrial reality and not the planning strategy used, it acts as a proxy for the current planning strategy
The \emph{Shortest} and \emph{Average} heuristics represent conventional strategies commonly applied in practice.
The \emph{Slope Scaling} and \emph{IP-based Local Search} represents state of the art heuristics from the SND literature.
They were adapted to our context, with detailed explanations in Appendix~\ref{appendix:benchmark_adaptation}.
For completeness, we also include the constructive heuristic and the lower bound rounding heuristic.
All benchmarks were implemented with identical code optimizations whenever applicable.  

Because the model introduced in this paper is new to the literature, no directly comparable algorithm is available.
On large-scale instances, the lower bound rounding heuristic serves as a natural mathematical programming baseline, consistent with rounding approaches used in flow problems \citep{jarrah2009, lamothe2021, lienkamp2024}.

\paragraph{Implementation}
We developed all algorithms in Julia \citep{julia} and released them as a public package, \texttt{ShipperTransportationPlanning.jl}.
To handle large-scale instances efficiently, we relied on several external libraries: \texttt{Graphs.jl} \citep{graphs} for custom graph structures and Dijkstra’s algorithm, \texttt{OhMyThreads.jl} \citep{ohmythreads} for parallelization, and \texttt{JuMP.jl} \citep{jump} with Gurobi for MILP solving. 
All experiments were conducted on a machine equipped with an Intel Core i9 processor (2.20 GHz) and 64 GB RAM.  

\subsection{Performance analysis} \label{subsec:perf_analysis}
\begin{figure}[!b]
\vspace{-0.25cm}
    \begin{subfigure}
        {0.34\textwidth}
        \centering
        \includegraphics[width=0.95\textwidth]{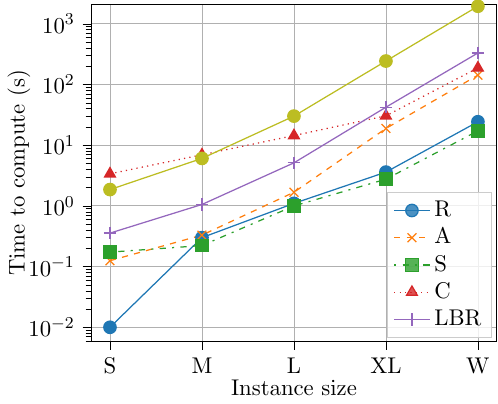}
        \caption{Time}
    \end{subfigure}
    \begin{subfigure}
        {0.33\textwidth}
        \centering
        \includegraphics[width=0.95\textwidth]{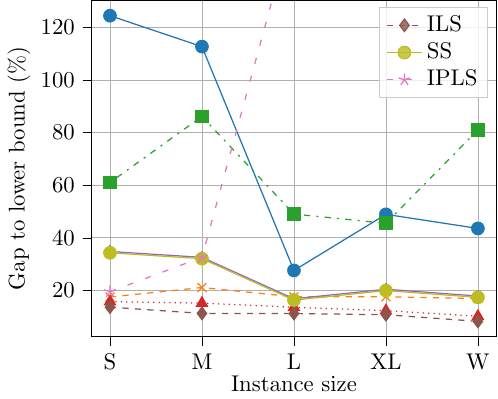}
        \caption{Gap}
    \end{subfigure}
    \begin{subfigure}
        {0.32\textwidth}
        \centering
        \includegraphics[width=0.95\textwidth]{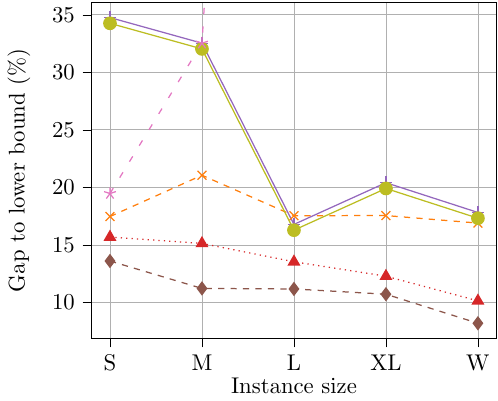}
        \caption{Gap (Zoom)}
    \end{subfigure}
    \caption{Results on scalability test instances}
    \label{fig:scale-results}
\end{figure}

\begin{figure}[!b]
    \begin{subfigure}
        {0.34\textwidth}
        \centering
        \includegraphics[width=0.95\textwidth]{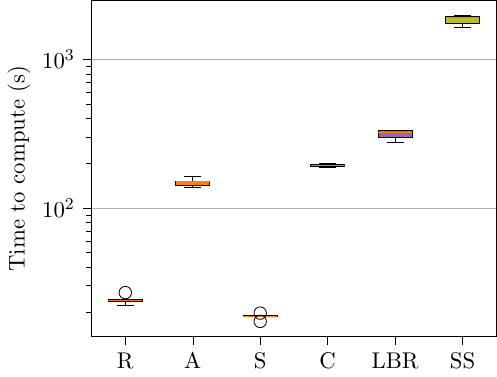}
        \caption{Time}
    \end{subfigure}
    \begin{subfigure}
        {0.325\textwidth}
        \centering
        \includegraphics[width=0.95\textwidth]{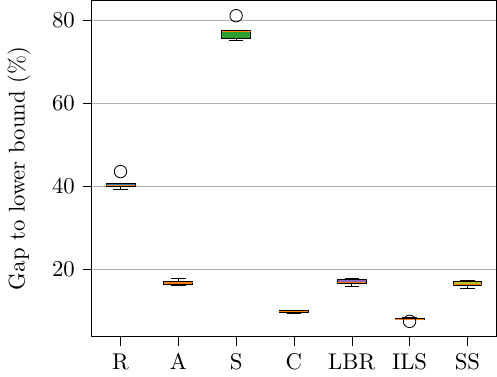}
        \caption{Gap}
    \end{subfigure}
    \begin{subfigure}
        {0.325\textwidth}
        \centering
        \includegraphics[width=0.95\textwidth]{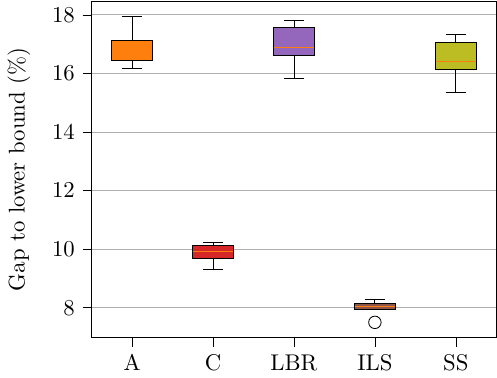}
        \caption{Gap (Zoom)}
    \end{subfigure}
    \caption{Results on performance test instances. 
    }
    \label{fig:perf-results}
\end{figure}

We evaluate how the proposed ILS algorithm improves over Renault’s proxy planning strategy and standard heuristic benchmarks, focusing on scalability and solution quality. 
Figures~\ref{fig:scale-results} and~\ref{fig:perf-results} compare ILS against all benchmarks on the two instance sets, reporting runtime and gap to the best available lower bound. 

Across all tests, ILS yields the highest solution quality, clearly outperforming other heuristics—especially on large instances—though at higher computational cost.
Benchmarks from the SND literature show limited effectiveness in this context: slope scaling yields only modest gains relative to lower‑bound rounding despite being slower, and IP‑based local search struggles to scale because solving exact subproblems is computationally demanding, reaching time or memory limits from the medium instance onward.

\begin{result}
    On the World instances, the Iterated Local Search performs best, with an average gap to the lower bound of 7.9\%.
    Further, it reveals a 23.2\% potential for Renault’s operational plans cost reduction. 
\end{result}

\begin{result}
    On the World instances, the Constructive heuristic performs second best, with an average gap to the lower bound of 9.9\%.
    It is 30\% faster than the lower bound rounding heuristic due to better parallelization. 
\end{result}

In practice, the choice between methods depends on context. 
For strategic planning in automotive inbound supply chains, where each percentage point represents substantial savings ($\sim$ 4 M\texteuro), ILS is preferred. 
For faster evaluations, such as tender analyses, the constructive heuristic is more appropriate. 
Simpler baselines like Shortest and Average remain informative but consistently underperform due to their inability to capture consolidation effects, confirming the need for tailored approaches such as ILS or the constructive heuristic.

\begin{result}
    Consolidation is necessary as illustrated by the poor performance of the Shortest heuristic.
    Average performs substantially better but is still 8.2\% worse than the ILS.
\end{result}  

\paragraph{Hyper-parameter tuning and Ablation study}

The hyper‑parameter study (in Appendix~\ref{appendix:performance_extensions}) shows that the local search consistently plateaus after about 45 minutes, and the full ILS converges after roughly six hours, which we adopt as the time limit for all experiments. 
Perturbations are kept lightweight ($\sim$2 minutes each) and constrained by variable and path thresholds to maintain tractability on large instances.
The ablation study (in Appendix~\ref{appendix:performance_extensions}) further confirms that the performance of the ILS does not rely on a single operator but on the complementary interaction of its components.

\subsection{Managerial and Practical Analyses} \label{subsec:ext_analysis}

In addition to validating the performance of our algorithm, we conduct further analyses to derive insights that are directly relevant for practitioners.
These analyses focus on three key aspects of transportation planning: the role of consolidation, transport regularity and outsourcing decisions.
Each of these dimensions reflects a fundamental design choice in Renault’s inbound supply chain and influences both computational tractability and managerial decision-making.

\paragraph{Consolidation model and Lower Bounds.}

We evaluate the impact of different relaxed consolidation models in Table~\ref{tab:bin_pack_consolidation}, which reports three lower bounds on instance Medium together with their rounded feasible solutions.
\begin{table}[!b]
\vspace{-0.75cm}
\centering
\caption{Comparison of lower bounds and associated rounding solutions on instance \textit{Medium}.} \label{tab:bin_pack_consolidation}
\renewcommand{\arraystretch}{1.3} 
\resizebox{\textwidth}{!}{
\begin{tabular}{>{\raggedright\arraybackslash}p{3.0cm} | >{\raggedright\arraybackslash}p{2.25cm} >{\raggedright\arraybackslash}p{2.25cm} >
{\raggedright\arraybackslash}p{2.25cm} >{\raggedright\arraybackslash}p{1.25cm} >{\raggedright\arraybackslash}p{2.5cm} >{\raggedright\arraybackslash}p{2.75cm} }
\toprule
Consolidation Type & Bound to bound gap & Sol. to bound gap & Sol. to sol. gap & Time & Sol. sensitivity to uncertainty & Sol. to sol. correspondence \\
\midrule
Linear Relaxation & -13.52\% & 64.50\% & 47.91\% & 1.0s & $\sim 3.7\%$ & $\sim 53\%$ \\
Mixed-Giant Cont. & -0.73\% & 32.54\% & 19.17\% & 1.1s & $\sim 4.1\%$ & $\sim 73\%$ \\
Full-Giant Cont. & 0.0\% & 35.31\% & 21.66\% & 3600s$^\star$ & $\sim 5.2\%$ & $\sim 69\%$ \\
Bin-Packing & - & 7.9\% & - & 1440s$^\star$ & $\sim 9.2\%$ & $100\%$ \\
\bottomrule
\end{tabular}
}
\note{The full-giant container bound computation was stopped after 1 hour.
It also provides a worst solution than the mixed-giant because it tends to fill the giant-container approximation, leading to excess bin with a real bin-packing.}
\end{table}
The full giant-container bound is the tightest, but the mixed relaxation offers the best trade-off, achieving comparable tightness much faster and producing the best rounded solution. 
Nevertheless, rounded solutions remain far from optimal. 
Beyond solution quality, the additional analyses reveal that consolidation fidelity also strongly affects network structure: over 30\% of giant-container paths differ from the bin-packing solution, with even lower correspondence for the linear relaxation. 
Common paths are mainly direct, high-volume, short-distance flows with limited consolidation potential. 
Although explicit bin-packing is more sensitive to uncertainty, this effect remains smaller than the structural differences between models, consistent with \citet{wang2019stochastic}.
These results confirm that bin-packing is both a computational bottleneck and a key determinant of solution quality, and that relaxations are useful for exploration (as in our ILS) but insufficient as standalone solution methods. 

\textbf{Managerial insight:} Accurate modeling of consolidation is essential. 
Approximating bin‑packing through aggregated capacities, as is common in the literature, can significantly distort solution quality. 
Even a one‑dimensional bin‑packing formulation, as used here, captures consolidation effects meaningfully while remaining computationally tractable.

\begin{figure}[!b]
    \begin{subfigure}{0.33\textwidth}
        \centering
        \includegraphics[width=0.95\textwidth]{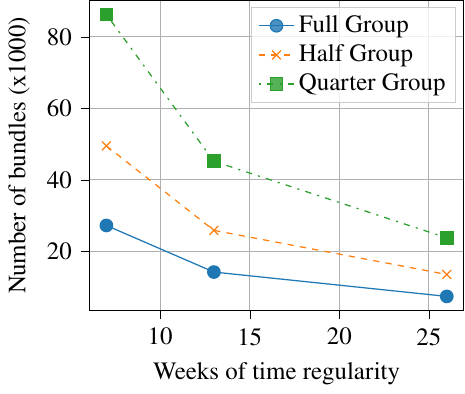}
        \caption{Instance size}
    \end{subfigure}
    \begin{subfigure}{0.33\textwidth}
        \centering
        \includegraphics[width=0.95\textwidth]{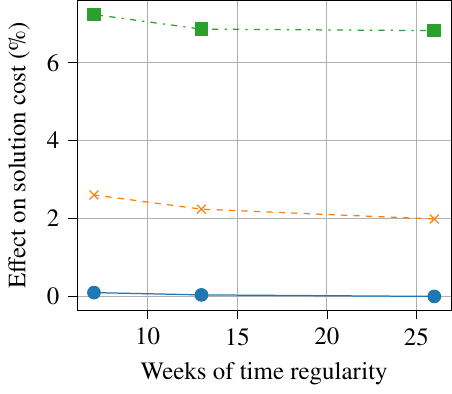}
        \caption{Constructive}
    \end{subfigure}
    \begin{subfigure}{0.33\textwidth}
        \centering
        \includegraphics[width=0.95\textwidth]{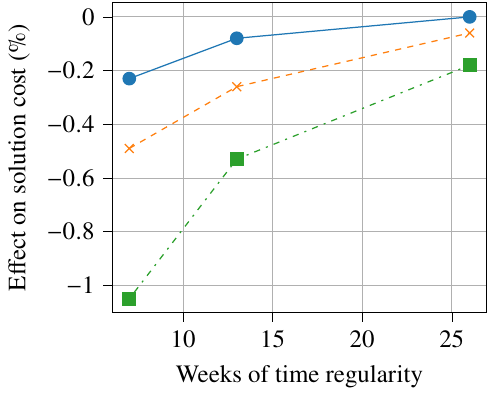}
        \caption{ILS}
    \end{subfigure}
    \caption{Impact of time regularity on the World instance.}
    \label{fig:regularity-results}
\end{figure}

\paragraph{Transport regularity.}
Figure~\ref{fig:regularity-results} shows the impact of time regularity on the World instance for both the constructive and the ILS heuristic.  
Lower regularity levels rapidly increase instance size and computational complexity. 
The constructive heuristic is highly sensitive to this effect, with solution quality deteriorating by up to 7\%. 
In contrast, the ILS benefits from the additional flexibility, achieving cost reductions of up to 1.05\%. 
This flexibility has only a marginal effect on the lower bound, which decreases by at most 0.25\%, consistent with the near‑linear nature of the relaxation.

\textbf{Managerial insight:} Regularity is valuable because it simplifies execution and reduces operational uncertainty. 
Our results confirm this: highly regular plans are easier to compute and remain cost‑competitive. 
While additional flexibility can provide small gains, the managerial trade‑off generally favors regular transportation plans, especially in high‑frequency inbound logistics.

\begin{figure}[!t]
    \begin{subfigure}
        {0.475\textwidth}
        \centering
        \includegraphics[width=0.9\textwidth, height=0.2\textheight]{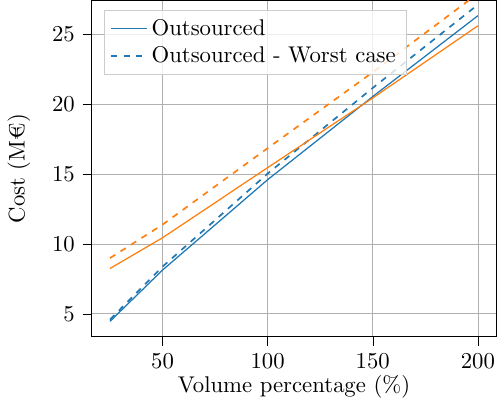}
        \caption{Medium}
    \end{subfigure}
    \begin{subfigure}
        {0.475\textwidth}
        \centering
        \includegraphics[width=0.9\textwidth, height=0.2\textheight]{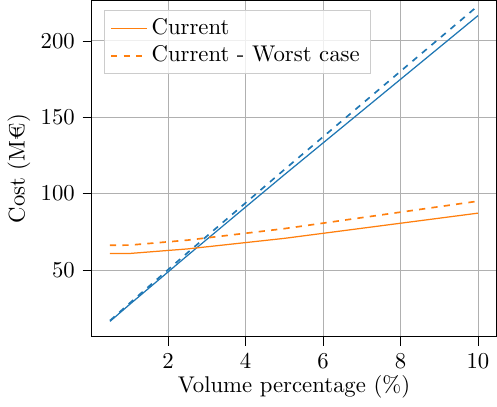}
        \caption{World}
    \end{subfigure}
    \caption{Impact of outsourcing on Medium and World instances.}
    \label{fig:outsource-results}
    \vspace{-0.25cm}
    \note{The volume of instance Medium is artificially expanded by duplicating commodities an adequate portion of the commodities.}
\end{figure}

\paragraph{Transport outsourcing.}
Figure~\ref{fig:outsource-results} compares the cost of a fully outsourced network to Renault’s current network across different volume configurations for the Medium and World instances. 
The outsourced baseline is constructed by replacing all non‑outsourced arcs with outsourcing costs estimated from average market rates. 
The results show that full outsourcing is attractive only at low volumes; as volumes increase, in‑house planning quickly becomes more cost‑effective. 
Larger networks reach this tipping point earlier, as they offer greater consolidation potential and stronger economies of scale.
In light of uncertainty, third-party logistics providers benefit from additional flexibility through capacity pooling. 
The worst-case scenarios for each strategy indicates that the insourcing strategy is more sensitive to uncertainty than outsourcing.

\textbf{Managerial insight:} The decision between outsourcing and internal planning depends critically on shipment volumes, network scale and sensitivity to uncertainty.
For smaller networks or low-volume settings, outsourcing may be cost-effective.
In contrast, for global, high-volume supply chains such as Renault’s, insourcing can yield substantial savings and greater strategic control.

\paragraph{Value of integrated planning and Solving to optimality}

Additional analyses reported in Appendix~\ref{appendix:managerial_extensions} further highlight the practical relevance of our approach. 
A decomposition experiment—separating intercontinental (strategic) and intra‑continental (tactical) planning—shows that such a two‑stage strategy consistently yields solutions about 2.5\% more expensive than the integrated model, confirming that consolidation, long‑haul routing, and timing decisions must be optimized jointly. 
We also examine the solvability of the proposed formulation with a state‑of‑the‑art MILP solver: while tiny instances can be solved to optimality, the Small instance already exceeds computational limits, underscoring the need for dedicated large‑scale heuristics such as the ILS developed in this paper.

\paragraph{Summary}
In summary, these analyses highlight four central lessons for both researchers and practitioners.  
First, accurate modeling of consolidation is indispensable.
Simplified representations 
yield computational speed, but they introduce distortions that can misguide planning decisions. 
Second, regularity emerges as a double-edged design lever.
Highly regular plans reduce problem size, enabling faster optimization, and align with planners’ preference for stable operations.
Third, outsourcing proves to be context-dependent.
Although outsourcing is attractive in small networks, low-volume flows or high uncertainty, large-scale inbound logistics—as faced by automotive inbound supply chains—benefit considerably from in-house planning and execution.
Fourth, integrated planning across long-haul and short-haul flows is essential: decomposing the strategic problem by flow scope overlooks the interactions between consolidation choices and timing decisions, and leads to different consolidation corridors.
Taken together, these insights demonstrate that our methodological advances are not purely of academic interest: they directly inform the design of shipper transportation networks in practice.

\section{Conclusion}
\label{sec:Conclusion}

This paper investigated the Shipper Transportation Planning Problem (STPP), a large-scale strategic design challenge motivated by Renault’s global inbound supply chain. 
We formalized a model that integrates three critical pillars of industrial logistics: explicit discrete consolidation, time-expanded routing, and operational regularity. 
To solve this problem at an unprecedented scale, we developed a tailored Iterated Local Search (ILS) metaheuristic. 
By combining large-neighborhood search with MILP-based perturbations and leveraging bundle-specific decompositions, the proposed framework provides high-quality feasible solutions alongside tractable lower bounds.

Our computational study on real industrial data demonstrates that the ILS significantly outperforms both legacy planning methods and adapted state-of-the-art heuristics from the service network design literature. On the largest ``World'' instances, which involve more than 700,000 commodities and 1.2 million arcs, the ILS identifies a potential 23.2\% reduction in transportation costs.
The figure compares Renault’s reconstructed historical plan with the optimized plan under identical assumptions, representing a theoretical rather than realized saving. At Renault’s scale, this potential reduction amounts to several hundred million euros annually.
Most notably, the framework has been successfully deployed in production at Renault, where it currently supports weekly strategic decisions. 
Forecast uncertainty, operational constraints, and gradual implementation limit the savings realized in practice. Renault’s internal evaluations nevertheless estimate annual benefits of approximately €20 million relative to its previous planning process, primarily from the highest-value transportation corridors, which are reviewed, negotiated with carriers, and progressively implemented.
Renault considers the deployment a major success and is progressively expanding the operational levers aligned with our model.
Beyond its practical impact, this work offers significant managerial and scientific insights:
\begin{description}
    \item[Consolidation Fidelity:] We demonstrate that moving beyond standard continuous capacity relaxations to explicit 1D bin-packing is a critical step forward; failing to capture the discrete nature of transportation units leads to substantial cost distortions.
    \item[The Price of Regularity:] We characterize the structural trade-offs between plan stability and cost, proving that highly regular plans satisfy planner preferences for stability with only marginal sacrifices in global efficiency.
    \item[Strategic Insourcing:] Our analysis identifies a volume-dependent threshold for outsourcing, suggesting that for global, high-volume shippers, in-house network design consistently outperforms third-party alternatives.
\end{description}

Future research could extend the STPP in several promising directions. 
Integrating multi-dimensional packing constraints (e.g., 3D stacking and axle weights) and the uncertainty inherent to forecasted data would further reduce the gap between strategic planning and dock-level execution.
A further direction is a hierarchical use of the framework, in which the corridors selected at the strategic level are passed to a tactical layer re-optimized at a finer temporal resolution with explicit dwelling, following the extensions of Appendix B.
Additionally, exploring how structural network properties---such as topology and spatial demand patterns---influence the performance of these large-scale systems would complement our practitioner-oriented findings.
In conclusion, this work provides a scalable algorithmic foundation and the first global-scale computational evidence for the shipper-side design of modern inbound supply chains.

\bibliographystyle{plainnat}
\bibliography{biblio}

\newpage

\appendix

\normalsize

These appendices provide additional material that complements the main text. 
Appendix~\ref{appendix:case_study} expands on the empirical case study by detailing the structure of Renault’s network and flow characteristics. 
Appendix~\ref{appendix:model_extensions} presents model extensions relevant for broader applications, while Appendix~\ref{appendix:lower_bound_proof} provides the formal proof of Proposition 2. 
Appendix~\ref{appendix:computational_efficiency} summarizes key implementation choices enabling computational scalability. 
Appendix~\ref{appendix:benchmark_adaptation} documents the adaptations required to apply benchmark SND heuristics to our setting. 
Appendix~\ref{appendix:performance_extensions} reports the full hyper‑parameter tuning and ablation study underlying the ILS configuration used in the experiments. 
Finally, Appendix~\ref{appendix:managerial_extensions} contains extended managerial and practical analyses—including integrated vs. decomposed planning and MILP‑solvability results—that further illustrate the practical implications of our approach.

\section{Detailed Case Study}\label{appendix:case_study}

\subsection{Supply Network}

Beyond the summary provided in the main paper, the full network maps in Figure~\ref{fig:sites-world-1} also show how Renault’s inbound flows span multiple continents and must be synchronized across a long‑haul platform core and a vast, sparsely connected supplier periphery. 
At this scale, lead times vary widely—from a few days for intra‑regional movements to several weeks for intercontinental shipments—making the system highly sensitive to external disruptions such as port congestion or geopolitical events. 
The detailed leg breakdown in Figure~\ref{fig:legs-world-1} further distinguishes supplier–site, supplier–platform, and platform–platform connections, highlighting the operational implications of each: tactical, short‑horizon management for supplier–site flows; strategic routing and frequency decisions for platform–platform links; and early capacity reservation for deep‑sea legs. 
These maps also clarify that, in the current project scope, each leg operates under a single available transport mode (trucks for inland flows and vessels for overseas movement), reflecting Renault’s data constraints, while still allowing multimodal extensions as outlined in Appendix B.

\subsection{Commodities and Flows}

\begin{figure}[!t]
    \begin{subfigure}{0.45\textwidth}
        \centering
        \includegraphics[width=0.95\textwidth]{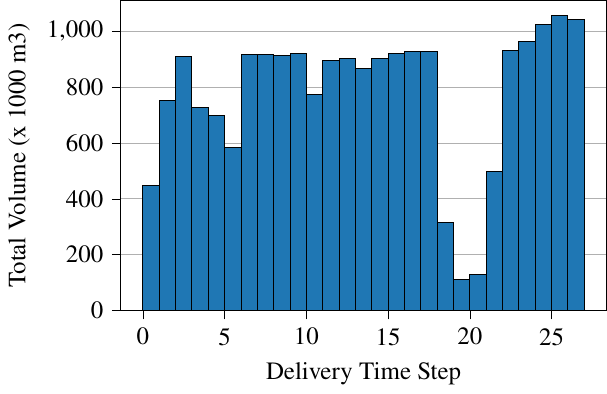}
        \caption{Volume per Week}
    \end{subfigure}
    \begin{subfigure}{0.45\textwidth}
        \centering
        \includegraphics[width=0.95\textwidth]{tikz_compiled/histogram_volume_com_size.pdf}
        \caption{Volume per Size}
    \end{subfigure}
    \caption{Commodities volume distributions}
    \label{fig:time-size-distribution}
\end{figure}

\begin{figure}[!t]
    \begin{subfigure}{0.45\textwidth}
        \centering
        \includegraphics[width=0.95\textwidth]{tikz_compiled/order_frequency_bundles2.pdf}
        \caption{Orders}
    \end{subfigure}
    \begin{subfigure}{0.45\textwidth}
        \centering
        \includegraphics[width=0.95\textwidth]{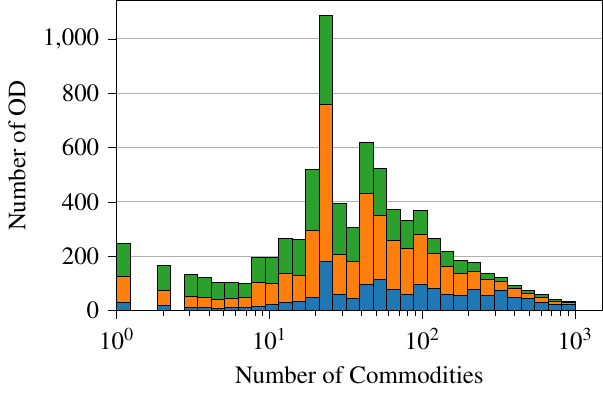}
        \caption{Commodities}
    \end{subfigure}
    \caption{Distribution of Orders and Commodities inside Origin-Destinations (OD)}
    \label{fig:order-com-inside-bundle}
\end{figure}

Renault’s inbound supply chain supports a highly complex and large-scale operation.
Over a six-month horizon, the network handles more than 40,000 distinct car parts, divided into 700,000 commodities to deliver, corresponding to approximately 9,000,000 packages and 20,000,000 m$^{3}$ of volume.
Figure~\ref{fig:time-size-distribution} illustrates the distribution of this total volume over time and across shipment sizes. 
The volume of inbound flows is relatively stable over time, without pronounced peaks or downswings.
This temporal regularity supports the implementation of structured mid-term planning processes.
For example, it enables the use of cyclic transport plans, fixed service frequencies, and pre-booked capacities for recurrent flows.
Such stability reduces the need for reactive short-term adjustments and improves predictability across the network.
However, the transportation plan must remain adaptable, as adjustments may occur frequently due to changes in supplier configurations or demand patterns.
Figure~\ref{fig:order-com-inside-bundle} illustrates the frequency of order and the number of commodities for origin-destination pairs, differentiating between short and long distances.
In this context, the sheer number of part types introduces substantial planning complexity.
The high granularity of the commodity mix leads to fragmented demand patterns, with each part potentially exhibiting different frequencies, sizes, and origin–destination characteristics.
This heterogeneity increases the dimensionality of planning problems and necessitates sophisticated aggregation, routing, and consolidation strategies to efficiently organize transportation.

Most individual shipments fall within a moderate size range of 1 to 4 cubic meters, as shown in figure (3b).
While these volumes are manageable, they typically do not fill a transport service on their own.
Consequently, high utilization of containers or trucks depends on consolidating flows across multiple dimensions—such as combining different part types, suppliers, or destinations within the same shipment.
This introduces additional decision layers involving timing, compatibility constraints, and routing coordination, especially for multi-stop or multi-leg routes.
The current organization of flows results in nearly 500,000 transport services—trucks or shipping containers—being used every six months.
Their spatial distribution is shown in Figure~\ref{fig:flows-world-1}.
Interestingly, only a small fraction of all potential legs in the network are actively used.
This is a direct consequence of the underlying consolidation strategy: flows are channeled through a limited number of high-volume corridors, where cost efficiency can be maximized.
In modeling terms, this indicates that leg activation is not static but endogenous—dependent on shipment densities, consolidation opportunities, and transport cost structures.
While most of the physical volume moves within Europe, the cost distribution exhibits a contrasting pattern.

\begin{figure}[!t]
    \begin{subfigure}{0.45\textwidth}
        \centering
        \includegraphics[width=0.95\textwidth]{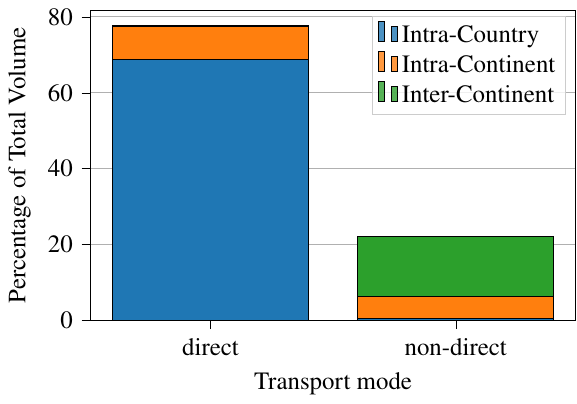}
        \caption{Volume repartition}
    \end{subfigure}
    \begin{subfigure}{0.45\textwidth}
        \centering
        \includegraphics[width=0.95\textwidth]{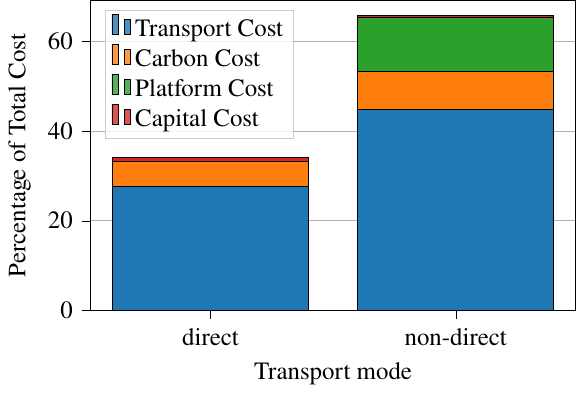}
        \caption{Cost repartition}
    \end{subfigure}
    \caption{Repartition of volumes and cost between direct and non-direct shipments}
    \label{fig:volume_cost_direct_distrib}
\end{figure}

Figure~\ref{fig:volume_cost_direct_distrib} illustrates the volume and cost repartition between flows shipped on direct legs and others.
The majority of transportation costs stem from intercontinental flows, which are fewer in number but far more expensive per service.
This cost–volume asymmetry has important implications for planning.
Intercontinental shipments require long lead times, early capacity reservation, and coordination with deep-sea or multimodal transport services.
In contrast, intra-European flows dominate operational complexity due to their higher frequency, network density, and time sensitivity.
Intercontinental flows are typically routed through designated global platforms that serve as strategic hubs for mode transitions and consolidation.
These hubs enable coordination across continents and help align long-haul shipments with downstream distribution schedules in Europe.
Their positioning and function are critical to synchronizing the global and regional components of the supply chain and mitigating the risk of disruptions or bottlenecks.

\section{Model Extensions and Discussions}\label{appendix:model_extensions}

As discussed in Section 3, the Shipper Transportation Planning model accommodates a wide range of extensions.
Those are described in the following paragraphs.

\paragraph{Multimodal Transport}

\begin{figure}[!b]
    \centering
    \includegraphics[width=0.75\linewidth]{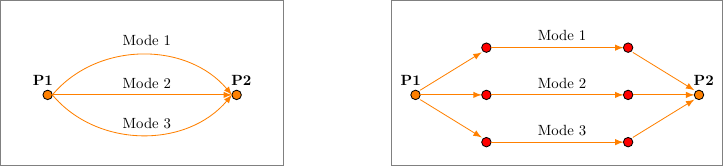}
    \caption{Example of multi-modal representation}
       \label{fig:multimodal}
\end{figure}

Few modes of transportation are used in Renault's problem: trucks for
inland transportation and boats for oversea transportation. Other modes could
be used, such as trains, barge or planes, to carry parts. To take this into
account, one could allow multiple arcs between two nodes, allowing for each arc to model one mode of transportation with different characteristics.
Another option is to add add dummy nodes
and arcs into the network. 
Figure~\ref{fig:multimodal} illustrates those two options.

\paragraph{Vehicle tours}

\begin{figure}[!t]
    \centering
    \includegraphics[width=0.75\linewidth]{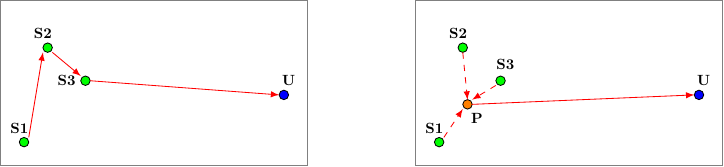}
    \caption{Example of vehicle tour approximation}
       \label{fig:vehicletour}
\end{figure}

Explicitly modeling vehicle tours, or milk-runs in the supply jargon,
between suppliers and plants is not included in the scope of this work. 
They can however be approximated without the need to
change our model by adding dummy platforms between proximate suppliers or plants.
Collecting or delivery arcs from this platform would be free while inter-platform arcs would cost the full tour.
Figure~\ref{fig:vehicletour} illustrates this option.

\paragraph{Inventory}

\begin{figure}[!b]
    \centering
    \includegraphics[width=0.75\linewidth]{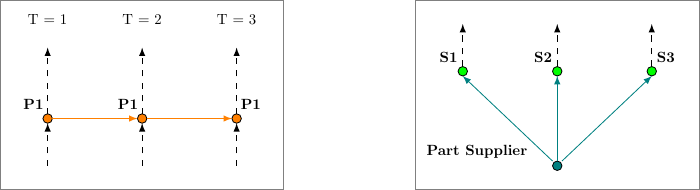}
    \caption{Example of inventory (left) and part sourcing (right) representation}
    \label{fig:inventory_part_sourcing}
\end{figure}

One can also notice that no mention of inventory is being made
throughout the model. This stems in our case from the coarse time step taken,
making all inventory happens inside time steps, which make them invisible.
Finer temporal resolution enables explicit inventory modeling by adding intra-node arcs between successive time steps, called \emph{inventory arcs}. 
An adequate \emph{inventory cost} would be applied on those arcs.
Figure~\ref{fig:inventory_part_sourcing} illustrates this.

\paragraph{Part Sourcing}

As previously mentioned, in our context, the same part can be provided by
several suppliers to several units but they are not considered substitutes
and are treated mathematically as different commodities. Integrating
part sourcing into this model can be done by adding a dummy node for
each part, like a super-supplier $S^{\mathrm{part}}$, connected to all eligible
suppliers of this part. 
Figure\ref{fig:inventory_part_sourcing} illustrates this. 
The flow
constraints on $\mathcal{G}$ would then need to be explicitly stated
, as
followed:
\[
    \sum_{a \in \delta^+(\nu)}f_{a}^{m}- \sum_{a \in \delta^-(\nu)}f_{a}^{m}
    = e_{\nu}^{m} \; \; \forall m \in M \text{, }\nu \in \mathcal{V}\cup
    S^{m}
\]
\[
    f_{\pi_\mathcal{G}(\alpha,o)}^{m}\leq C^{m}x_{\alpha}^{b}\quad \forall
    b \in B \text{, }o \in O_{b}\text{, }m \in M_{o}
\]

\paragraph{Return logistics}

By considering packaging used in plants as commodities to be due to
suppliers, adding a dummy supplier for each plant and a dummy plant for each
supplier, this model can tackle simultaneously forward and return
logistics of parts.

\paragraph{Temporal regularity extension}

One can extend the regularity constraint to the temporal dimension, enforcing that all cargoes belonging to the same supplier–unit pair be shipped at a predefined frequency.
Let $n$ denote the desired shipment periodicity for a given supplier-unit pair. 
Temporal regularity requires that all orders belonging to this pair depart only on dates congruent to a fixed subset of $\{1,\ldots,n\}$ modulo $n$. 
A natural way to enforce this structure is to refine the bundle construction. Specifically, we split each bundle associated with this supplier–unit pair into $n$ sub-bundles, corresponding to the $n$ congruence classes modulo $n$. The $k$-th sub-bundle contains all orders whose delivery date satisfies
\(
t \equiv k \pmod{n}.
\)
By design, such a partition ensures that each sub-bundle contains only orders compatible with one admissible departure pattern.
To prevent shipments from leaving on forbidden departure dates, we modify the time-expanded transportation graph: any arc corresponding to a disallowed departure time for the sub-bundle in question is assigned an infinite cost (or removed from the feasible arc set). 
Since admissible arcs are determined by projecting travel times into the weekly time-expanded network, this modification preserves the temporal feasibility of all remaining paths.

\paragraph{Limited number of services}

If a hard upper bound on the number of services available on an arc is required, it can be incorporated in the heuristic framework with minimal changes. 
This limit can be handled by detecting arcs where bin-packing would require more than $K_a^{\max}$ services: when this occurs, the algorithm simply marks the corresponding timed arc in the travel-time or time-space graph as infeasible (equivalently, assigns it an infinite cost).
A limitation arise when this limited number of services renders the feasible path set empty for a bundle.
Addressing this situation would require the use of a strong penalization strategy instead of a hard infeasibility, allowing controlled violations of feasibility in the constructive heuristic and subsequent recovery within the ILS.
This mechanism automatically propagates through the constructive heuristic, local search, and perturbation phases, preventing any path that violates hard service capacities.
In practice, this modification integrates seamlessly into the current cost-evaluation routines, and the bin-packing subproblems can be adapted to stop early and signal infeasibility whenever the service limit is exceeded.

\paragraph{Service Frequency Requirements}

Transportation contracts may impose minimum service frequencies or regular departure schedules. 
Such requirements arise when a carrier commits to operating a transportation service a prescribed number of times during a planning period, independently of the total volume transported. 
For example, a weekly demand transported under a contract requiring three departures per week must be distributed across three services. 
This requirement can be accommodated through a straightforward extension of the consolidation model: instead of applying the bin-packing procedure to the aggregate volume $\sum_{m \in M} f_a^m \ell_m$ transported on arc $a$, the flow is first partitioned into $n$ service streams $f_{a,1}, \dots, f_{a,n}$ corresponding to the required departures, and bin-packing is then performed independently on each stream. 
The transportation cost therefore becomes 
\[ 
    C_a^{\mathrm{freq}}(f_a) = \sum_{i=1}^{n} C_a^{\mathrm{con}}(q_{a,i}) 
\] 
which accounts for the vehicle requirements associated with each scheduled departure while preserving the same underlying consolidation mechanism.

\section{Proof of Proposition 2}\label{appendix:lower_bound_proof}

We consider the linear relaxation of problem (9). 
For $a \in \mathcal{A}^{\mathrm{con}}$ and $k \in K_a$, as $\tau_a^k \in \mathbb{R}_+$, we have: 
$$ L_a\tau_{a}^{k}= \sum_{m\in M}y_{ak}^{m}\ell_{m}$$
For $a \in \mathcal{A}^{\mathrm{con}}$ and $m \in M$, bin-packing constraints imposes: 
$$ \sum_{k \in K_a} y_{ak}^{m} = f_{a}^{m}$$
We have therefore for all $a \in \mathcal{A}^{\mathrm{con}}$: 
$$ \sum_{k \in K_a}\tau_{a}^{k} = \sum_{m \in M}f_{a}^{m}\frac{\ell_{m}}{L_a}$$

In this relaxation, the number of transport units used by commodities on each arc equals the total commodity volume on the arc divided by the unit capacity. 
These equations transforms packing constraints into ``liquid container'' volume
cost, accounted for in the linearly relaxed network cost function in problem (16).
\begin{equation*}
    \sum_{a \in \mathcal{A}} \sum_{m \in M}(c_{am}^{\mathrm{com}}+ \frac{\ell_{m}}{L_{a}}c
    _{a}^{\mathrm{con}}) f_{a}^{m} + \sum_{p \in \mathcal{P}}c_{p}
    ^{\mathrm{over}}z_{p}\label{obj:linearNetCost}
\end{equation*}

Without platform overloading constraints (7), this linear relaxation of bin-packing constraints allows decoupling between commodities, which is equivalent in our case to a decoupling between bundles. 
Our shortest path constraint is written as a flow constraint on the graph
$\mathscr{G}$. 
Because the polytope of flows is perfect, the value we obtain is the linear relaxation bound.
\Halmos
    
\section{Computational efficiency considerations}\label{appendix:computational_efficiency}

The practical scalability of our algorithm relies on five key implementation features, in addition to the efficiency of the underlying mathematical formulation. 
We summarize these components below.
\begin{description}
    \item[\textbf{Bin‑packing heuristics.}] To avoid solving exact bin‑packing problems during each cost evaluation, we rely on a fast First‑Fit Decreasing (FFD) heuristic, which computes packings within milliseconds and has proven sufficiently accurate for our purposes.
    \item[\textbf{Problem reduction.}] We exploit a refined version of the lower‑bound rounding heuristic to prune bundles that will never be routed through the common network. 
    Specifically, any bundle whose direct‑shipping heuristic cost is lower than the best lower‑bound cost involving the common network is removed from consideration.
    \item[\textbf{Pre‑computations.}] Two types of pre‑computations are used to reduce online processing time: (i) For each bundle, we pre‑identify the subset of feasible arcs, avoiding unnecessary cost evaluations on the entire network, and (ii) Empty bin‑packing configurations are pre‑computed and reused, which significantly accelerates cost evaluation on arcs through which no flow is routed.
    \item[\textbf{Memory management.}] We use shared memory structures tailored to repeated bin‑packing evaluations. 
    In particular, commodity sets are converted into integer vectors, and a single pre‑allocated buffer is reused across evaluations, avoiding costly memory allocation during local‑search iterations.
    \item[\textbf{Parallelism.}] Parallel computation is applied at two levels: bin‑packing evaluations within the RePack neighborhood and arc‑cost computations within the Insert operator. 
    Because these tasks require memory isolation, we rely on thread‑specific memory buffers (via Julia “Channels”), achieving both concurrency and safe memory usage.
\end{description}

\section{Adaptation of SND benchmarks}\label{appendix:benchmark_adaptation}

Classical heuristics developed for Service Network Design (SND) cannot be applied directly to our setting, as the standard SND modeling framework does not fully align with the structure of the Shipper Transportation Planning Problem. 
In particular, SND formulations typically assume a flat network representation, predefined service paths, and linear or fixed‑charge cost structures, whereas our problem features explicit bin‑packing consolidation, bundle‑level regularity constraints, and a massive path space that cannot be enumerated. 
Nonetheless, to enable a meaningful comparison with the SND literature, we adapted two widely used heuristic approaches: the slope‑scaling algorithm of \citet{jarrah2009} and the IP‑based local search of \citet{lindsey2016improved}.

\paragraph{Slope‑Scaling Heuristic}

The slope‑scaling algorithm is based on iteratively linearizing and rescaling cost coefficients to approximate the nonlinear components of the objective.
In our problem, nonlinearity arises from the bin‑packing structure, which determines the number of services required on each arc. 
To adapt the method, we replace the bin‑packing cost on each consolidated arc aaa with a rescaled volumetric cost following the update rule inspired by Jarrah et al.: $$ c_a^{(k)} \;=\; \frac{\mathrm{BinPacking}(a)}{\mathrm{TotalVolume}(a)} $$ where $\mathrm{BinPacking}(a)$ denotes the number of services required on arc $a$ under the current packing solution and $\mathrm{TotalVolume}(a)$ is the aggregated commodity volume routed through $a$. 
At each iteration, these scaled costs $c_a^{(k)}$ define a linearized problem that is solved to generate updated flows, which in turn determine new bin‑packing estimates. 
This adaptation preserves the spirit of the original algorithm while making it compatible with our consolidation framework.

\paragraph{IP‑Based Local Search}

The IP‑based local search method of Lindsey, Erera, and Savelsbergh relies on solving a sequence of exact mixed‑integer subproblems that reoptimize parts of the solution. 
In our adaptation, the MILP formulation used for each subproblem is directly inherited from the full model presented in Section~\ref{sec:PbStatement}, with the restriction that only a selected subset of bundles is reoptimized at each iteration. 
This mirrors the logic of the perturbation subproblems in our ILS, with one key difference: in the IP‑based local search, the bin‑packing constraints are modeled exactly, rather than through the giant‑container relaxations used within the perturbation phase of the ILS.
As a result, each subproblem provides a high‑fidelity local improvement, at the cost of significantly greater computational effort.

\section{ILS hyper-parameter tuning and ablation study}\label{appendix:performance_extensions}

We propose here the detail of the ILS hyper-parameter tuning to help readers better understand how those choices were made and the ablation study to help understand how each element contributes to the global performance.

\begin{figure}[!b]
    \begin{subfigure}{0.33\textwidth}
        \centering
        \includegraphics[width=0.95\textwidth]{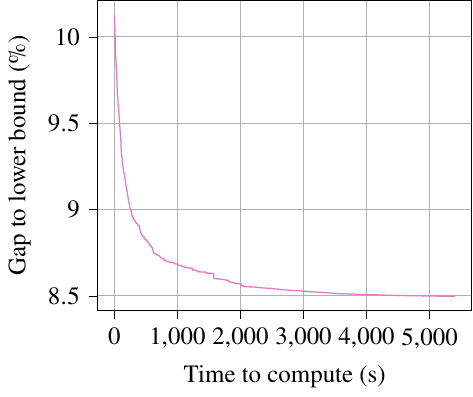}
        \caption{Local search profile}
    \end{subfigure}
    \begin{subfigure}{0.33\textwidth}
        \centering
        \includegraphics[width=0.95\textwidth]{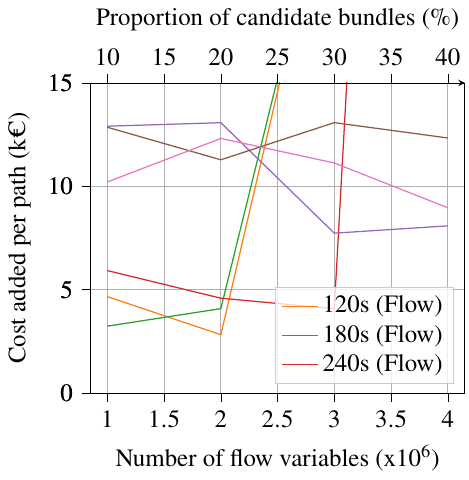}
        \caption{Perturbation impact}
    \end{subfigure}
    \begin{subfigure}{0.33\textwidth}
        \centering
        \includegraphics[width=0.95\textwidth]{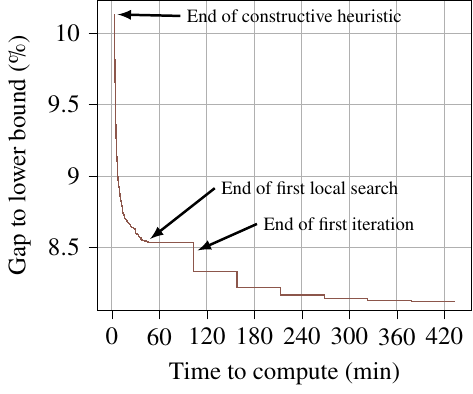}
        \caption{ILS profile}
    \end{subfigure}
    \caption{Influence of parameters on local search and perturbation and the resulting ILS.}
    \label{fig:hyperparameters-profiles}
\end{figure}

\paragraph{ILS hyper-parameters and time profiles}\label{subsec:hyperparameters}
Several hyperparameters guide the execution of the ILS algorithm, as described in Section~\ref{sec:ILS}. 
The most important are the time limits for the local search and the perturbation phase. 
We calibrate perturbations further using two parameters: the maximum number of MILP variables and the minimum number of candidate paths. 
Figure~\ref{fig:hyperparameters-profiles} illustrates the influence of these parameters on both components and the resulting profile of the ILS obtained.  
Figure~\ref{fig:hyperparameters-profiles} (a) shows that the local search start to converge after approximately 45 minutes. 
We therefore set the time limit for each local search in the ILS to 45 minutes.  
To balance the scale of perturbations with the degradation they introduce, we allocate two minutes for each perturbation. 
We restrict flow-based perturbations to 2{,}000{,}000 variables and require path-based perturbations to affect at least 30\% of the bundles. 
We repeat perturbations until they modify 15\% of the paths or increase the solution cost by more than 2\%.
Figure~\ref{fig:hyperparameters-profiles} (c) reports the time profile of the ILS with these parameter settings on the \emph{World} instance. As can be seen, the ILS converges after six hours of computation, corresponding to roughly six full iterations, and continues to improve only marginally thereafter.
Accordingly, we adopt a six-hour time limit for all subsequent experiments in the computational study.

\paragraph{Ablation study.}
To assess the individual contribution of each component of the proposed ILS, we performed an ablation study in which operators were selectively deactivated. 
We considered four ablation possibilities: (i) \emph{no re-insert}, (ii) \emph{no consolidate-and-refine}, (iii) \emph{no single-plant perturbation} and (iv) \emph{no attract-reduce perturbation}.
We do not consider the \emph{no re-pack} possibility because this component is the basis of every other. 
Not using it is akin to not using the local search nor the perturbation scheme.  
We consider the full combination if these ablated configurations and each variant was run using identical stopping criteria.

\begin{figure}[!b]
    \begin{subfigure}
        {0.34\textwidth}
        \centering
        \includegraphics[width=0.95\textwidth]{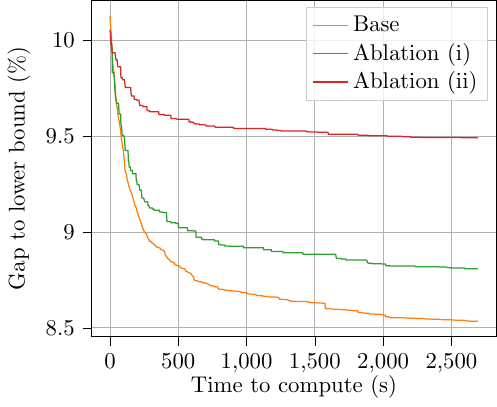}
        \caption{Local Search}
    \end{subfigure}
    \begin{subfigure}
        {0.33\textwidth}
        \centering
        \includegraphics[width=0.95\textwidth]{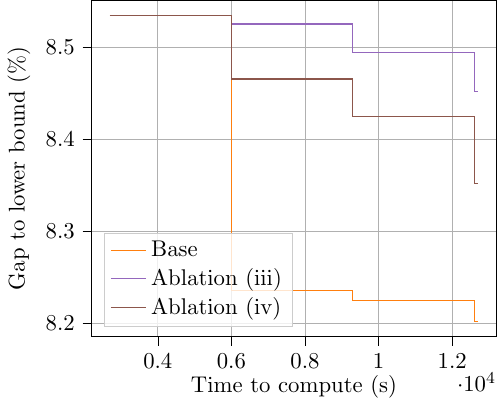}
        \caption{Perturbations}
    \end{subfigure}
    \begin{subfigure}
        {0.32\textwidth}
        \centering
        \includegraphics[width=0.95\textwidth]{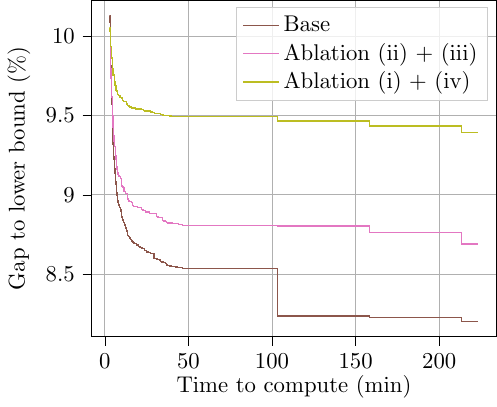}
        \caption{ILS}
    \end{subfigure}
    \caption{Results of the ablation study tests}
    \label{fig:ablation-study}
\end{figure}

Figure~\ref{fig:ablation-study} compares the time profile  of each variant, all starting from the same greedy solution.

The results show that every building block of the ILS contributes positively to the overall performance. 
In the local search, removing the Re‑Insert neighborhood degrades cost-reduction effectiveness by 62\% (the final cost is 1.2\% worse), while removing the Consolidate neighborhood leads to a 25\% effectiveness degradation (with a 0.5\% increase in cost).
Both neighborhoods prove essential on large instances, improving consolidation quality and reducing the number of required services.

The perturbation layer exhibits the least significant effect: restricting the algorithm to a single perturbation mechanism leads to a less efficient exploration mechanism, hindering the overall performance by 0.25\% for the \emph{Single-Plant} perturbation and 0.15\% for the \emph{Attract-Reduce} perturbation. 

It is clear that the overall effects of the different operators stacks.
Overall, the ablation study confirms that the strength of the proposed ILS derives not from a single dominant component but from the complementary interaction of multiple neighborhoods and perturbation strategies.

\section{Extended Managerial and Practical Analyses}\label{appendix:managerial_extensions}

\paragraph{Value of integrated planning.}
Renault’s inbound logistics network encompasses transportation legs with heterogeneous operational characteristics, ranging from short-haul, intra‑regional supplier flows to long‑haul intercontinental services. 
Usually, such heterogeneity often motivates a separation of planning responsibilities across strategic and tactical levels. 
To quantify the implications of this common practice, we complement our computational study with a decomposition-based approach in which the proposed ILS heuristic is applied sequentially.
First to the intercontinental (strategic) portion of the network, and subsequently to the intra-continental (tactical) legs, treating the strategic outcome as fixed input for the second stage.
Table~\ref{tab:decomposition} reports the solution quality and computational effort associated with both approaches. 
\begin{table}[!b]
\centering
\caption{Comparison of solution quality with and without decomposition on performance instances.}\label{tab:decomposition}
\renewcommand{\arraystretch}{1.3} 
\begin{tabular}{>{\raggedright\arraybackslash}p{3.0cm} | >{\raggedright\arraybackslash}p{3.25cm} >{\raggedright\arraybackslash}p{2.75cm} >
{\raggedright\arraybackslash}p{2.75cm} >{\raggedright\arraybackslash}p{1.25cm} }
\toprule
Setting & Average gap & Minimum gap & Maximum gap & Time \\
\midrule
Integrated & 7.9\% & 7.5\% & 8.3\% & 6h \\
Decomposed & 10.4\% & 9.7\% & 10.8\% & 1h30 \\
\bottomrule
\end{tabular}
\end{table}
The results indicate that although the decomposed strategy reduces computation time, it consistently produces inferior solutions: across benchmark instances, the integrated model yields solutions that are, on average, 2.5\% less costly than the best decomposed alternative. 
In effect, the decomposition removes the performance gains brought by the ILS scheme relative to a greedy baseline.

\textbf{Managerial insight:}  Despite its operational convenience, the decomposition strategy performs worst when applied at scale, where consolidation decisions and departure schedules require joint coordination. 
The performance gap arises from the interplay between consolidation, long‑haul routing, and timing decisions interdependencies that are only fully captured under integrated planning. 
Feedback from Renault further confirms the importance of maintaining this coupling in practice.

\paragraph{Solving to optimality}
It is interesting to explore to what extent the proposed problem can be solved to optimality using a state-of-the-art MILP solver. 
To do so, we introduce two additional benchmark instances, referred to as Tiny (T) and Extra‑Small (XS), whose structural characteristics are summarized in Table~\ref{tab:tiny-instances}.
\begin{table}[!t]
\centering
\caption{Description of instances \emph{Tiny} and \emph{Extra-Small}.}\label{tab:tiny-instances}
\renewcommand{\arraystretch}{1.3} 
\begin{tabular}{>{\raggedright\arraybackslash}p{1.5cm} | >{\raggedright\arraybackslash}p{1.5cm} >{\raggedright\arraybackslash}p{1.5cm} | >
{\raggedright\arraybackslash}p{1.5cm} >{\raggedright\arraybackslash}p{1.75cm} | >{\raggedright\arraybackslash}p{1.5cm} >{\raggedright\arraybackslash}p{1.5cm} >{\raggedright\arraybackslash}p{1.5cm} >{\raggedright\arraybackslash}p{1.5cm}}
\toprule
 & \textbf{$\mathscr{G}$ nodes} & \textbf{$\mathscr{G}$ arcs} & \textbf{$\mathcal{G}$ nodes} & \textbf{$\mathcal{G}$ arcs} & \textbf{$B$} & \textbf{$O$} & \textbf{$M$} & \textbf{$K$} \\
\midrule
T & 111 & 475 & 156 & 812 & 21 & 54 & 950 & 9 \\
\midrule
XS & 157 & 773 & 576 & 3 270 & 23 & 123 & 1547 & 98 \\
\bottomrule
\end{tabular}
\end{table}

Table~\ref{tab:optimality-results} reports the computational performance of Gurobi on these instances, as well as on the original Small instance. 
\begin{table}[!t]
\centering
\caption{Description of Problem~\eqref{pb:forwardPath} resolution with Gurobi.}\label{tab:optimality-results}
\renewcommand{\arraystretch}{1.3} 
\begin{tabular}{>{\raggedright\arraybackslash}p{1.5cm} | >{\raggedright\arraybackslash}p{2.5cm} >{\raggedright\arraybackslash}p{2.5cm} | >
{\raggedright\arraybackslash}p{3.5cm} | >
{\raggedright\arraybackslash}p{3.5cm}}
\toprule
Instance & \textbf{\# variables} & \textbf{\# constraints} & \textbf{Time to optimal sol.} & \textbf{ILS gap to bound} \\
\midrule
T & 116 800 & 26 800 & 38 s. & 0\% \\
\midrule
XS & 5 970 300 & 218 600 & 5760 s. & 0 \% \\
\midrule
S & 248 430 500 & 23 250 100 & - & - \\
\bottomrule
\end{tabular}
\note{We consider here the time to which the solver found the optimal solution but has not yet proved its optimality.}
\end{table}
The results highlight the sharp growth in model size as instance dimensions increase. 
While instances T and XS can be solved to optimality, within seconds for T and within hours for XS, the Small instance already exceeds practical solvability limits.
The results further indicate that, for these smaller instances, the greedy heuristic attains solutions that match the optimal bound. 
However, from the Small instance onward, direct MILP solution becomes prohibitive, reinforcing the need for dedicated heuristics capable of handling realistic problem sizes.

\paragraph{Sensitivity to aggregate platform-capacity contracts.} 

Constraint~(8) models the contractual volume threshold specified in platform agreements.
While constraint~(7) limits the volume that can be processed during a given period, constraint~(8) captures the cumulative volume above which overload tariffs apply over the contract horizon. 
To evaluate its impact, we re-optimized all benchmark instances after removing constraint~(8). 
The resulting transportation plans remained almost unchanged: more than $99.7\%$ of bundle paths were preserved, corresponding to roughly 20 modified bundles among 7,000, and the objective value changed by only $0.6\%$. Furthermore, when the aggregate threshold is evaluated \emph{a posteriori} on the unconstrained solutions, only one platform exceeds its contracted volume, by $1.9\%$ of its capacity. For comparison, Renault's current transportation plan exceeds the threshold on four platforms, representing $5.4\%$ of total contracted capacity. These results indicate that aggregate platform-capacity contracts have little influence on routing and consolidation decisions. This outcome is consistent with the cost structure of the case study, where transportation procurement accounts for approximately $80\%$ of total costs and therefore dominates corridor-selection decisions. Nevertheless, the aggregate constraint improves the accuracy of platform-cost estimation and allows the model to correctly account for overload tariffs negotiated in practice. Its value therefore lies primarily in procurement and contract evaluation rather than in shaping the transportation network itself.
 
\end{document}